\documentclass[11pt]{article}
\usepackage[T1]{fontenc}

\usepackage{amsfonts}
\usepackage{amssymb}
\usepackage{latexsym}
\usepackage{color}
\usepackage{graphicx}
\usepackage{eurosym}
\usepackage[latin1]{inputenc}
\usepackage{amsmath}
\usepackage{natbib}
\usepackage{mathtools}

\usepackage{setspace}
\usepackage{fancyhdr}\usepackage{ifthen}
\makeatletter
\def\@biblabel#1{}
\renewcommand\@cite[2]{{#1\if@tempswa,\nolinebreak[3] #2\fi}}
\makeatother

\newtheorem{thm}{Theorem}[section] 
\newtheorem{ass}[thm]{Assumption} 
\newtheorem{prop}[thm]{Proposition} 
\newtheorem{coro}[thm]{Corollary} 
\newtheorem{lem}[thm]{Lemma} 
\newtheorem{rem}[thm]{Remark} 
\newtheorem{exa}[thm]{Example} 
\numberwithin{equation}{section}

\newcommand{\bfc}{\mbox{$\mbox{\boldmath $c$}$}}
\newcommand{\bfd}{\mbox{$\mbox{\boldmath $d$}$}}
\newcommand{\bfe}{\mbox{$\mbox{\boldmath $e$}$}} 
\newcommand{\bff}{\mbox{$\mbox{\boldmath $f$}$}} 
\newcommand{\bfg}{\mbox{$\mbox{\boldmath $g$}$}} 

\newcommand{\bfp}{\mbox{$\mbox{\boldmath $p$}$}} 
\newcommand{\bfs}{\mbox{$\mbox{\boldmath $s$}$}} 
\newcommand{\bft}{\mbox{$\mbox{\boldmath $t$}$}} 
\newcommand{\bfu}{\mbox{$\mbox{\boldmath $u$}$}} 
\newcommand{\bfw}{\mbox{$\mbox{\boldmath $w$}$}} 
\newcommand{\bfx}{\mbox{$\mbox{\boldmath $x$}$}} 
\newcommand{\bfy}{\mbox{$\mbox{\boldmath $y$}$}} 
\newcommand{\bfz}{\mbox{$\mbox{\boldmath $z$}$}} 

\newcommand{\sbff}{\mbox{\boldmath \scriptsize $f$}} 
\newcommand{\sbfg}{\mbox{\boldmath \scriptsize $g$}} 

\newcommand{\sbfp}{\mbox{\boldmath \scriptsize $p$}} 

\newcommand{\sbxi}{\mbox{\boldmath \scriptsize $\xi$}} 

\newcommand{\bfF}{\mbox{$\mbox{\boldmath $F$}$}} 
\newcommand{\bfG}{\mbox{$\mbox{\boldmath $G$}$}}
\newcommand{\bfH}{\mbox{$\mbox{\boldmath $H$}$}}

\newcommand{\bfU}{\mbox{$\mbox{\boldmath $U$}$}} 

\newcommand{\bfW}{\mbox{$\mbox{\boldmath $W$}$}}
\newcommand{\bfX}{\mbox{$\mbox{\boldmath $X$}$}} 

\newcommand{\bdelta}{\mbox{$\mbox{\boldmath $\delta$}$}} 
\newcommand{\bfeta}{\mbox{$\mbox{\boldmath $\eta$}$}} 

\newcommand{\bmu}{\mbox{$\mbox{\boldmath $\mu$}$}} 
\newcommand{\bxi}{\mbox{$\mbox{\boldmath $\xi$}$}} 
\newcommand{\bomega}{\mbox{$\mbox{\boldmath $\omega$}$}}

\newcommand{\bXi}{\mbox{$\mbox{\boldmath $\Xi$}$}} 

\newcommand{\bfzero}{\mbox{$\mbox{\boldmath $0$}$}} 
\newcommand{\bfone}{\mbox{$\mbox{\boldmath $1$}$}}

\newcommand{\bigzero}{\scalebox{1.0}{\mbox{\large $0$}}}

\date{}

\begin{document}

{\textbf{\Large \rule{0cm}{1cm} \hspace{-0.225cm}THE GRANGER\textendash JOHANSEN REPRESENTATION \\
\rule{0cm}{0.6cm}THEOREM FOR INTEGRATED TIME SERIES\\
\rule{0cm}{0.6cm}ON BANACH SPACE}}

\textbf{Phil Howlett}\footnote{\textbf{Corresponding Author: Phil Howlett}, Scheduling and Control Group (SCG), Centre for Industrial and Applied Mathematics (CIAM), UniSA STEM, University of South Australia. email:~phil.howlett@unisa.edu.au, url:~http://orcid.org/0000-0003-2382-8137.}, \textbf{Brendan K Beare}\footnote{\textbf{Brendan K Beare}, School of Economics, University of Sydney. email:~brendan.beare@sydney.edu.au, url:~https://orcid.org/0000 -0001-9146-131X.}, \textbf{Massimo Franchi}\footnote{\textbf{Massimo Franchi}, Sapienza University of Rome, Rome, Italy. email:~massimo.franchi@uniroma1.it, url:~https://orcid.org/0000-0002-3745-2233.}, \textbf{John Boland}\footnote{\textbf{John Boland}, Centre for Industrial and Applied Mathematics, UniSA STEM, University of South Australia. email:~john.boland@unisa.edu.au. url:~https://orcid.org/0000-0003-1132-7589.}, \\
\textbf{Konstantin Avrachenkov}\footnote{\textbf{Konstantin Avrachenkov}, INRIA, Sophia Antipolis, France. email:~k.avrachenkov@sophia.inria.fr, url:~https://orcid.org/0000-0002-8124-8272}. \\

\textbf{Abstract:} We prove an extended Granger\textendash Johansen representation theorem (GJRT) for finite or infinite order integrated autoregressive time series on Banach space. We assume only that the resolvent of the autoregressive polynomial for the series is analytic on and inside the unit circle except for an isolated singularity at unity. If the singularity is a pole of finite order the time series is integrated of the same order. If the singularity is an essential singularity the time series is integrated of order infinity. When there is no deterministic forcing the value of the series at each time is the sum of an almost surely convergent stochastic trend, a deterministic term depending on the initial conditions and a finite sum of embedded white noise terms in the prior observations. This is the extended GJRT.  In each case the original series is the sum of two separate autoregressive time series on complementary subspaces\textemdash a singular component which is integrated of the same order as the original series and a regular component which is not integrated.  The extended GJRT applies to all integrated autoregressive processes irrespective of the spatial dimension, the number of stochastic trends and cointegrating relations in the system, and the order of integration.

\textbf{Keywords:} $AR$ time series, resolvent operators, GJRT, cointegrated series. \\
\rule{0cm}{0.6cm}\textbf{JEL:} C32. \\
\rule{0cm}{0.6cm}\textbf{MSC\! [2020]:} 62M10, 91B84, 47A11, 47A55. \\

\setcounter{footnote}{0}

\vspace{-0.5cm}

\section{Introduction}
\label{s:int}

We wish to establish an extended form of the Granger\textendash Johansen representation theorem (GJRT) for vector autoregressive ($AR$) time series taking values in a Banach space. These series are often referred to as functional time series because, for instance, each real-valued integrable function on the unit interval of the real line can be represented by an infinite-dimensional vector whose components are the coefficients of the Fourier sine series. In this particular case the sum of the magnitudes of the vector components will be finite. Thus we can study the space of integrable functions on the unit interval by considering the vector space of all infinite-dimensional vectors where the sum of the magnitudes of the components is finite. Consequently the study of vector $AR$ processes on Banach space can be regarded as a generalized form of the study of $AR$ time series for functions.

An important early contribution to the analysis of functional time series is \citet{bos1}, where a theoretical treatment of linear processes in Banach and Hilbert spaces is developed. In particular the derivation of laws of large numbers and central limit theorems allows estimation and inference for infinite-dimensional stationary $AR$ models. These models also enable direct representation of the dynamics of infinite-dimensional objects, such as a time series of continuous functions on a compact spatial domain, and they allow greater generality for modelling of conditional means and variances \citep{horv1}.

Economic applications of functional time series include studies on the term structure of interest rates \citep{kar1}, intraday volatility \citep{horm1, gab1} and the human influence on climate \citep{cha1}. Additional applications can be found in the recent monograph by \cite{kok2}.

One might expect the analysis of functional time series to begin by testing for stationarity. In this regard \cite{horv2}, \cite{kok1} and \cite{kok2} suggest extended forms of the Kwiatkowski-Phillips-Schmidt-Shin (KPSS) test of stationarity proposed originally for single variable series by \cite{kwi1} while \cite{aue1} devise a test of stationarity in the frequency domain. However the literature on this topic is somewhat sparse. An alternative approach by \cite{cha2} proposed a test of non-stationarity based on the generalized eigenvalues for the covariance operator of the observed values and on the long-run covariance operator for the associated first differences. If the dimension of the non-stationary subspace is finite then the dimension is called the number of common stochastic trends. In a subsequent working paper \cite{hu1} consider an infinite-dimensional $AR(1)$ process with a compact operator and show that the common trends representation comprises a finite number of stochastic trends each integrated of order one and an infinite-dimensional cointegrating space. Hu and Park also propose an estimator for the functional autoregressive operator. More recently \cite{nie1} suggest a statistical procedure to determine the dimension of the nonstationary subspace of cointegrated functional time series taking values in a Hilbert space of square-integrable functions defined on a compact interval. As in \cite{cha2} they assume that the cointegrated series has a finite-dimensional nonstationary subspace while the stationary subspace is infinite-dimensional. That is, there are infinitely many linearly independent cointegrating relations. The test is applied to several empirical examples: age-specific U.S. employment rates, Australian temperature curves, and Ontario electricity demand. For more details of functional time series we refer to extended reviews in \cite{bea2}, \cite{fra5} and \cite{nie1}.

Our primary concern here is to review the development of the GJRT for $AR$ processes. According to \cite{han1} the GJRT for a finite-dimensional cointegrated vector $AR(1)$ time series $\bfx \sim I(1)$ can be stated informally as follows.

\begin{quote}
\label{q:han}
{\em A cointegrated vector autoregressive process can be decomposed into four components: a random walk, a strictly stationary process, a deterministic part, and a term that depends on the initial conditions.}
\end{quote}

The original theorem proposed by \cite{eng1} was a representation theorem for finite-dimensional cointegrated vector processes $\bfx \sim I(1)$ that connected the moving average, autoregressive, and error correction forms. The focus for Johansen was to find a more explicit representation for finite-dimensional $AR(p)$ processes $\bfx = \{ \bfx(t)\}_{t \in {\mathbb Z}}$ satisfying
\begin{equation}
\label{e:(1.1)}
\bfx(t) = \Phi_1\bfx(t-1) + \cdots + \Phi_p \bfx(t-p) + \bxi(t)
\end{equation}
for all $t \in {\mathbb Z}$ where $\Phi_1,\ldots,\Phi_p \in {\mathbb R}^{n \times n}$ are matrix coefficients and $\bxi(t) \in {\mathbb R}^n$ is a strong white noise process. The specific focus on $AR(p)$ processes enabled Johansen to find necessary and sufficient conditions for $\bfx \sim I(1)$ \citep{joh1} and for $\bfx \sim I(2)$ \citep{joh2}. In each case the conditions depend only on the properties of the matrix coefficients $\Phi_1,\ldots,\Phi_p$.

\cite{bea1} used the model described in (\ref{e:(1.1)}) for a process taking values in a Hilbert space $H$. In this case $\bfx(t), \bxi(t) \in H$ for all $t \in {\mathbb Z}$ and $\Phi_1,\ldots,\Phi_p \in {\mathcal B}(H)$ are bounded linear operators. Beare et al.~defined an autoregressive polynomial $\Phi:{\mathbb C} \rightarrow {\mathcal B}(H)$ by the formula
\begin{equation}
\label{e:(1.2)}
\Phi(z) = I - \mbox{$\sum_{i=1}^p$} \Phi_i z^i
\end{equation}
for each $z \in {\mathbb C}$.  For $p=1$ they show that if the spaces $\Phi(1)(H)$ and $\Phi^{\, \prime}(1)\Phi(1)^{-1}( \{\bfzero \})$ are closed complementary subspaces in $H$ then the resolvent $\Phi(z)^{-1}$ has a simple pole at $z = 1$ and $\bfx \sim I(1)$. For $p > 1$ they reach the same conclusion using the additional assumption that $\Phi_1,\ldots,\Phi_p$ are compact. In a subsequent paper \cite{bea2} assume that $\Phi_1,\ldots,\Phi_p$ are compact and show that $\Phi(z)$ has a simple pole at $z=1$ if and only if $\Phi(1)(H)$ and $\Phi^{\, \prime}(1)\Phi(1)^{-1}( \{\bfzero \})$ are closed complementary subspaces in $H$.  They also show that if $\Phi(z)^{-1}$ is analytic for $z \in D_{1+ \epsilon}(0) \setminus \{1\} = \{z \in {\mathbb C } \mid |z| < 1 + \epsilon \} \setminus \{1\}$ for some $\epsilon > 0$ and has a pole of order two at $z=1$ with
\begin{equation}
\label{e:(1.3)}
\Phi(z)^{-1} = \Upsilon_{-2}/(z-1)^2 + \Upsilon_{-1}/(z-1) + \mbox{$\sum_{\ell \in {\mathbb N}-1}$} \Upsilon_{\ell}(z-1)^{\ell}
\end{equation}
for all $z \in D_{0,\, \epsilon}(1) = \{z \in {\mathbb C} \mid 0 < |z-1| < \epsilon \}$ then $\bfx \sim I(2)$. They use the Moore\textendash Penrose inverse to derive a formula for the coefficient $\Upsilon_{-2}$ and a more complicated formula for $\Upsilon_{-1}$. \cite{fra5} extend these results to show that if $\Phi(z)^{-1}$ has a pole of order $d$ at $z=1$ with a corresponding unit root of finite type then $\bfx \sim I(d)$. The assumption that $\Phi_1,\ldots,\Phi_p$ are compact and the alternative assumption that the unit root for $\Phi(z)$ is of finite type are each sufficient to ensure that the subspace $\Upsilon_{-1}(H) \subset H$ is finite-dimensional. Thus determination of the generalized eigenspace for $\Phi(1)$ and the associated Laurent series coefficients is essentially reduced to a matrix problem. The procedure proposed by \cite{fra5} is essentially equivalent to earlier work by \cite{avr2}, \cite{avr3}, \cite{fra1} and \cite{how1} on the inversion of matrix power series.  See also \cite{how2}.  The work by \cite{bea1}, \cite{bea2} and \cite{fra5} is restricted to Hilbert spaces where $\Upsilon_{-1}(H)$ is finite-dimensional.  

\cite{seo1} extends the results in \cite{bea2} to $AR$ processes $\bfx \sim I(1) \vee I(2)$ taking values in a Banach space $X$. The main results depend on two critical assumptions\textemdash for $\bfx \sim I(1) \vee I(2)$, that the subspaces $\Phi(1)(X)$ and $\Phi(1)^{-1}(\{\bfzero\})$ can each be complemented in the space $X$ and, for $\bfx \sim I(2)$, that certain additional key subspaces in $X$ can also be complemented. Complementation of closed subspaces is not guaranteed  in Banach space. Consequently \cite{seo1} leaves major theoretical questions unanswered. Seo considers only first and second order poles of the resolvent operator at $z = 1$ and does not consider isolated essential singularities.

It is well known that a necessary and sufficient condition for an $AR$ process $\bfx$ to be integrated of order $d < \infty$ is that $\Phi(z)^{-1}$ has a pole of order $d$ at $z=1$. In the early papers on cointegration a plethora of seemingly {\em ad hoc} conditions were used collectively to determine the order of integration. In some sense the collective nature of these conditions obscured the fact that the underlying necessary and sufficient condition was simply the order of the pole.

The aim of this paper is to obtain an extended form of the GJRT that is valid for all $AR$ processes on Banach space.  In this regard it is sufficient to establish the extended form for first order vector autoregressive $AR(1)$ processes\footnote{Any $AR(p)$ process on a space $X$ can be modelled as an $AR(1)$ process using augmented autoregressive coefficients on the product space $X^p$.} on Banach space. We will use recent research on the inversion of operator pencils \citep{alb2,alb5} that has not previously been used for the analysis of integrated time series.

By confining our attention to the naturally complemented subspaces defined by the spectral projections for the autoregressive polynomial at $z=1$, we remove the need for restrictive assumptions about complementation of certain subspaces or their finite dimensions as in \cite{bea1}, \cite{bea2}, \cite{fra5}, and \cite{seo1}.  The theory presented in this paper applies to all integrated autoregressive processes irrespective of the dimension of the space, the number of stochastic trends and cointegrating relations in the system, and the order of integration.  Hence the present results are a substantial extension of those reported in \cite{bea1}, \cite{bea2}, \cite{fra5}, and \cite{seo1}.  We provide several particular examples to illustrate the key structural issues including a composite graph showing the simulated evolution of an archetypal $I(\infty)$ process.  

The rest of the paper is organised as follows. Section~\ref{s:mr} outlines the key components of the new representation and summarizes the significance of our contribution, Section~\ref{s:rrmt} describes the mathematical tools used in the solution procedure and in Section~\ref{s:rtsm} these tools are applied to solution of the time series model.  Section~\ref{s:rtsm} also includes formal statements of the main results.  Section~\ref{s:sc} shows that existing $I(1)$ and $I(2)$ representations are special cases of the extended GJRT.  The general structures are illustrated in Section~\ref{ss:tepr} via three examples and in Section~\ref{s:aepg} the new methods are applied to solve a non-autonomous problem for dynamic control of electrical power generation in a wind farm. In Section~\ref{s:c} we draw some brief conclusions. In the Appendix our notation and terminology are listed and explained in Section~\ref{s:nt}, the proofs of all formally stated results are given in Section~\ref{s:pmr} and the relevant properties of the weighted Volterra (integral) operator used in the wind farm problem are justified in Section~\ref{s:wvo}.

\section{Main results}
\label{s:mr}

Our main contribution, presented formally as Proposition~\ref{p:1}, is an extended  GJRT on Banach space that applies to integrated $AR(1)$ series of any spatial dimension, any number of stochastic trends and cointegrating relations, and any order of integration. Let $A_0, A_1 \in {\mathcal B}(X,Y)$ and suppose that $\bfx = \{ \bfx(t)\}_{t \in {\mathbb Z}} \in X^{\mathbb Z}$ satisfies
\begin{equation}
\label{e:(2.1)}
A_0 \bfx(t) + A_1 \bfx(t-1) = \bxi(t)
\end{equation}
for all $t \in {\mathbb N}-1$ where $\bfx(-1) = \bfc$ and $\bxi = \{\bxi(t)\}_{t \in {\mathbb Z}} \in Y^{\mathbb Z}$ is a strong white noise. The autoregressive polynomial is $A(z) = A_0 + A_1z = B_0 + B_1(z-1) \in {\mathcal B}(X,Y)$ for all $z \in {\mathbb C}$ where $B_0 = A_0 + A_1 \in {\mathcal B}(X,Y)$ and $B_1 = A_1 \in {\mathcal B}(X,Y)$ and the resolvent is $R(z) = A(z)^{-1} \in {\mathcal B}(Y,X)$ for all $z \in {\mathbb C} \setminus \sigma$ where $\sigma$ is the spectral set for $A(z)$.

\begin{ass}
\label{a:1}
The resolvent $R(z) = A(z)^{-1}$ is analytic for $z \in [D_1(0)^a \setminus \{1\}] \cup\, D_{0,\delta}(1) = \{z \in {\mathbb C} \mid |z| \leq 1 \mbox{ or } 0 < |z-1| < \delta \} \setminus \{1\}$ for some $\delta > 0$ and has an isolated singularity at $z = 1$.
\end{ass}

\textbf{\textit{The extended representation:}} The series $\bfx = \{\bfx(t)\}_{t \in {\mathbb N}-1}$ is the sum of a stochastic trend
\begin{equation}
\label{e:(2.2)}
\mbox{$\sum_{k \in {\mathbb N}}$} (-1)^k T_{-k} \Delta^{-k} \bxi_+(t),\qquad T_{-k} = (-1)^{k-1}(T_{-1}B_0)^{k-1}T_{-1},
\end{equation}
where $\Delta$ is the difference operator, $\Delta^{-1}$ is the inverse difference or cumulation operator, $\bxi_+(t) = \bxi(t)$ for $t \in {\mathbb N}-1$ and $\bxi_+(t) = \bfzero$ otherwise, and $T_{-1}$ is the residue of $R(z)$ at $z=1$, a deterministic component
\begin{equation}
\label{e:(2.3)}
(-1)R_t B_1 \bfc,\qquad R_t = (-1)^t(A_0^{-1}A_1)^tA_0^{-1},
\end{equation}
and a finite weighted sum of strong white noise terms
\begin{equation}
\label{e:(2.4)}
\mbox{$\sum_{s=0}^t$} W_s\, \bxi(t-s),\qquad W_s = R_s + (I_X - T_{-1}B_0)^{-s-1}T_{-1}.
\end{equation}

\textbf{\textit{Additional observations:}} Our new results extend the GJRT in various ways.

\begin{enumerate}

\item If $R(z)$ has an isolated singularity at $z = 1$ then Lemma~\ref{lem:rcon} shows that $R(z)$ is analytic for $z \in [ D_1(0)^a \setminus \{1\}] \cup\, D_{0,\delta}(1)$ for some $\delta > 0$ if and only if there is some $\epsilon > 0$ such that $R(z)$ is analytic for $z \in D_{1+\epsilon}(0) \setminus \{1\}$.

\item The Laurent series $R(z) = \sum_{j \in {\mathbb Z}} T_j(z-1)^j$ for $z \in D_{0,\epsilon}(1)$ is found by solving the fundamental equations (\ref{e:(3.1)}) and (\ref{e:(3.2)}) for $\{T_j\}_{j \in {\mathbb Z}}$ subject to the magnitude constraints (\ref{e:(3.3)}). The resolvent $R(z)$ is analytic for $z \in D_{0,\epsilon}(1)$ if and only if the fundamental equations have a basic solution $\{T_{-1}, T_0 \}$.  The resolvent is completely determined by the basic solution. A general solution procedure for the fundamental equations is available in separable Banach space but simpler methods can always be used in finite dimensional problems.

\item  There is a closed formula $R(z) = [(z-1)I_X + T_{-1}B_0]^{-1}T_{-1} + [I_X + T_0B_1(z-1)]^{-1}T_0$ for the resolvent when $z \in D_{0, \epsilon}(1)$.  The singular part is $R_{\mbox{\scriptsize \rm sin}}(z) = [(z-1)I_X + T_{-1}B_0]^{-1}T_{-1} = \sum_{k \in {\mathbb N}} T_{-k}/(z-1)^k$ when $z \neq 1$ where $T_{-k} = (-1)^{k-1}(T_{-1}B_0)^{k-1}T_{-1}$.  The regular part is $R_{\mbox{\scriptsize \rm reg}}(z) = [I_X + T_0B_1(z-1)]^{-1}T_0 = \sum_{\ell \in {\mathbb N}-1} T_{\ell}(z-1)^{\ell}$ when $z \in D_{\epsilon}(1)$ where $T_{\ell} = (-1)^{\ell} (T_0B_1)^{\ell}T_0$.  These formul{\ae} are new to time series analysis. 

\item The coefficients $(-1)^kT_{-k}$ for the loadings on the stochastic trend components, $R_t$ for the deterministic components and $W_s$ for the weights of the embedded noise terms are easily calculated using (\ref{e:(2.2)}), (\ref{e:(2.3)}) and (\ref{e:(2.4)}).

\item The resolvent $R(z)$ may have a pole of order $d \in {\mathbb N}$ at $z=1$ with $\bfx \sim I(d)$ or an isolated essential singularity at $z=1$ with $\bfx \sim I(\infty)$. In the former case $T_{-1}B_0$ is nilpotent of degree $d$ and the stochastic trend reduces to a finite sum. In the latter case $T_{-1}B_0$ is quasi-nilpotent with $(T_{-1}B_0)^k \neq \bigzero$ for all $k \in {\mathbb N}$ but with $\Vert (T_{-1}B_0)^k \Vert^{1/k} \to 0$ as $k \to \infty$.

\item The operators $P = T_{-1}B_1 \in {\mathcal B}(X)$ and $Q = B_1T_{-1} \in {\mathcal B}(Y)$ are the spectral separation projections for $A(z)$ at $z=1$.  The direct sum decompositions $X = P(X) \oplus P^c(X)$ and $Y = Q(Y) \oplus Q^c(Y)$ separate the singular and regular parts of $R(z)$ with $R_{\mbox{\scriptsize \rm sin}}(z)A(z) = P$ and $A(z)R_{\mbox{\scriptsize \rm sin}}(z) = Q$ while $R_{\mbox{\scriptsize \rm reg}}(z)A(z) = P^c$ and $A(z)R_{\mbox{\scriptsize \rm reg}}(z) = Q^c$. It follows that $T_{-k} \in {\mathcal B}(Q(Y),P(X))$ for all $k \in {\mathbb N}$ and $T_{\ell} \in {\mathcal B}(Q^c(Y),P^c(X))$ for all $\ell \in {\mathbb N}-1$.

\item The stochastic trend lies entirely in $P(X)$ and depends only on previous values of the noise component in $Q(Y)$. The weighted sum of embedded noise terms lies entirely in $P^c(X)$ and depends only on previous values of the noise component in $Q^c(Y)$. The term $R_tB_1 \bfc$ can be written as the sum of $U_tB_1P\bfc \in P(X)$ and $W_t B_1 P^c \bfc \in P^c(X)$.

\item The $AR(1)$ time series $\bfx \in X^{\mathbb Z}$ can be expressed as the sum of two separate $AR(1)$ time series $\bfx_{\mbox{\scriptsize{sin}}} \in [P(X)]^{\mathbb Z}$ and $\bfx_{\mbox{\scriptsize{reg}}} \in [P^c(X)]^{\mathbb Z}$. The singular component $\bfx_{\mbox{\scriptsize{sin}}}$ is defined by $(QA_0P) \bfx_{\mbox{\scriptsize{sin}}}(t) + (QA_1P) \bfx_{\mbox{\scriptsize{sin}}}(t-1) = Q \bxi(t)$ for $t \in {\mathbb Z}$ with $\bfx_{\mbox{\scriptsize{sin}}}(-1) = P \bfc$ and the regular component $\bfx_{\mbox{\scriptsize{reg}}}$ is defined by $(Q^cA_0P^c) \bfx_{\mbox{\scriptsize{reg}}}(t) + (Q^cA_1P^c) \bfx_{\mbox{\scriptsize{reg}}}(t-1) = Q^c \bxi(t)$ for $t \in {\mathbb Z}$ with $\bfx_{\mbox{\scriptsize{reg}}}(-1) = P^c \bfc$.  The singular component $R_{\mbox{\scriptsize{sin}}}(z) = PR(z)Q$ is analytic for all $z \neq 1$ and the extended regular component $R_{\mbox {\scriptsize \rm reg}}(z) = P^cR(z)Q^c$ is analytic for $z \in D_{1+\epsilon}(0)$.  If $P \bfc = \bfzero$ then $\bfx_{\mbox{\scriptsize{sin}}} \sim I(d)$ is integrated of order $d$ for some $d \in {\mathbb N} \cup \{ \infty\}$.  If $P^c \bfc = \sum_{s \in {\mathbb N}-1} \Phi^s \bfeta(-1-s)$ where $A_{i,r} = Q^cA_iP^c \in {\mathcal B}(P^c(X), Q^c(Y))$ for each $i=0,1$ and $\Phi = (-1)A_{0,r}^{-1}A_{1,r} \in {\mathcal B}(P^c(X))$ then $\bfx_{\mbox{\scriptsize{reg}}} \sim I(0)$ is a standard linear process.   If both conditions hold\footnote{These sufficient conditions were not stated explicitly in the original paper: Phil Howlett, Brendan K. Beare, Massimo Franchi, John Boland and Konstantin Avrachenkov (2025), The Granger-Johansen representation theorem for integrated time series on Banach space,  \textit{Journal of Time Series Analysis}, \textbf{46}, 432\textendash 457, http://doi.org/ 10.1111/jtsa.12766.} then $\bfx = (\bfx_{\mbox{\scriptsize{sin}}} + \bfx_{\mbox{\scriptsize{reg}}}) \sim I(d)$ is also integrated of order $d$.

\item The key subspaces $P(X), P^c(X), Q(Y), Q^c(Y)$ need not be finite-dimensional and there are no restrictions on the operators $B_0, B_1 \in {\mathcal B}(X,Y)$. The subspaces $T_{-1}(Y) \subset P(X)$ and $[T_{-1}(Y)]^{\perp} \subset X^*$ are the attractor space and cointegrating space respectively.

\item If $\bfx \in {\mathbb R}^n$ then $R(z) = \mbox{adj}[B_0 + B_1(z-1)]/\det[B_0 + B_1(z-1)] \in {\mathbb R}^{n \times n}$ and the form of the Laurent series $R(z) = \sum_{j \in {\mathbb N -n}} T_j(z-1)^j$ is known with $T_{-k} = \bigzero$ for $k > n$.  The matrix formula gives no direct indication that $R(z)$ is completely determined by $\{T_{-1}, T_0\}$.

\item Our decision to restrict discussion to {\em strictly stationary} noise processes is one of convenience. Analogous results remain true if the noise process is {\em weakly stationary}. In Banach space a process $\bxi \in X^{\mathbb Z}$ is said to be weakly stationary \citep[Definition 2.4, p 65]{bos1} if ${\mathbb E}[ \|\bxi(t)\|^2] < \infty$ for all $t \in {\mathbb Z}$, ${\mathbb E}[\bxi(t)] = \bmu$ does not depend on $t \in {\mathbb Z}$ and
$$
{\mathbb E} \left[ \langle \bxi(s+t) - \bmu, \bff \rangle \overline{\langle \bxi(s) - \bmu, \bfg \rangle} \right] = c_{\sbff, \,\sbfg}(t)
$$
depends only on $t \in {\mathbb Z}$ for each fixed pair of linear functionals $\bff, \bfg \in X^*$ and all $s \in {\mathbb Z}$. 
\end{enumerate}

\section{Representation of the resolvent and the \boldmath ${\mathcal M}$-transform}
\label{s:rrmt}

We consider the $AR(1)$ process defined in (\ref{e:(2.1)}). The autoregressive polynomial is a bounded linear operator pencil $A(z) = A_0 + A_1z = B_0 + B_1(z-1) \in {\mathcal B}(X,Y)$ for all $z \in {\mathbb C}$.  \cite{alb2} showed that the resolvent $R(z) = A(z)^{-1} \in {\mathcal B}(Y,X)$ exists and can be expressed as a Laurent series $R(z) = \mbox{$\sum_{j \in {\mathbb Z}}$} T_j(z-1)^j$ for all $z \in D_{0,\,\epsilon}(1)$ if and only if the coefficients $\{T_j\}_{j \in {\mathbb Z}} \in {\mathcal B}(Y,X)^{\mathbb Z}$ satisfy a system of left and right fundamental equations
\begin{equation}
\label{e:(3.1)}
T_{j-1}B_1 + T_jB_0 = \left \{ \begin{array}{ll}
I_X & \mbox{if}\ j=0 \\
\bigzero_X & \mbox{if}\ j \in {\mathbb Z}, j \neq 0 \end{array} \right.
\end{equation}
and
\begin{equation}
\label{e:(3.2)}
B_1T_{j-1} + B_0T_j = \left \{ \begin{array}{ll}
I_Y & \mbox{if}\ j=0 \\
\bigzero_Y & \mbox{if}\ j \in {\mathbb Z}, j \neq 0 \end{array} \right.
\end{equation}
where the coefficients $\{T_j\}_{j \in {\mathbb Z}}$ also satisfy the magnitude constraints
\begin{equation}
\label{e:(3.3)}
\lim_{k \rightarrow \infty} \|T_{-k}\|^{1/k} = 0 \quad \mbox{and} \quad \lim_{\ell \rightarrow \infty} \|T_{\ell}\|^{1/\ell} \leq 1/\epsilon.
\end{equation}
If (\ref{e:(3.1)}), (\ref{e:(3.2)}) and (\ref{e:(3.3)}) are all satisfied then
\begin{equation}
\label{e:(3.4)}
T_{-k} = (-1)^{k-1}(T_{-1}B_0)^{k-1}T_{-1} \in {\mathcal B}(Y,X)
\end{equation}
for all $k \in {\mathbb N}$ and
\begin{equation}
\label{e:(3.5)}
T_{\ell} = (-1)^{\ell} (T_0B_1)^{\ell}T_0 \in {\mathcal B}(Y,X)
\end{equation}
for all $\ell \in {\mathbb N}-1$. Thus the series is completely determined by the basic solution $\{T_{-1}, T_0 \}$. Now we can write $R(z) = R_{\mbox{\scriptsize \rm sin}}(z) + R_{\mbox{\scriptsize \rm reg}}(z)$ where the singular part of the Laurent series
\begin{equation}
\label{e:(3.6)}
R_{\mbox{\scriptsize \rm sin}}(z) = [I_X(z-1) + T_{-1}B_0]^{-1}T_{-1} = \mbox{$\sum_{k \in {\mathbb N}}$} T_{-k}(z-1)^{-k}
\end{equation}
converges for all $z \neq 1$ and the regular part of the Laurent series
\begin{equation}
\label{e:(3.7)}
R_{\mbox{\scriptsize \rm reg}}(z) = [I_X + T_0B_1(z-1)]^{-1}T_0 = \mbox{$\sum_{\ell \in {\mathbb N}-1}$} T_{\ell}(z-1)^{\ell}
\end{equation}
converges for all $z \in D_{\epsilon}(1)$.  Formal statements of these results can be found in \cite{alb2} and \citet{alb5}.  The next result shows that the necessary and sufficient conditions for the above representation are satisfied if and only if $R(z)$ is analytic for $z \in [D_1(0)^a \setminus \{ 1\}] \cup D_{0, \delta}(1)$ for some $\delta > 0$ and has an isolated singularity at $z = 1$. 

\begin{lem}
\label{lem:rcon}
Assumption \ref{a:1} holds if and only if $R(z)$ is analytic on a set $D_{1 + \epsilon}(0) \setminus \{1\}$ for some $\epsilon > 0$ and has an isolated singularity at $z=1$. 
\end{lem}

Our preferred method of solution is the Maclaurin transform, which we shall refer to as an ${\mathcal M}$-transform.  The ${\mathcal M}$-transform of the stochastic time series $\bfu \in X^{\mathbb Z}$ with ${\mathbb E}[\|\bfu(\omega,t)\|] = \int_{\Omega} \| \bfu(\omega,t) \| d \mu(\omega) \leq c/r^t$ for some $c, r > 0$ and all $t \in {\mathbb N}-1$ is an almost surely convergent random power series denoted by ${\mathcal M}[\bfu](z) = \bfU(z):\Omega \rightarrow X$ and defined by
\begin{equation}
\label{e:(3.8)}
\bfU(\omega, z) = \mbox{$\sum_{t \in {\mathbb N}-1}$} \bfu(\omega, t) z^t
\end{equation}
for each $z \in D_r(0) \subset {\mathbb C}$ and all $\omega \in \Omega$.  We justify the almost sure convergence as follows.  Choose $z \in D_r(0)$.  Now ${\mathbb E}[ \| \bfU(\omega, z) \|] \leq \sum_{t \in {\mathbb N}-1} {\mathbb E}[ \| \bfu(\omega, t) \|] |z|^t \leq \sum_{t \in {\mathbb N}-1} c(|z|/r)^t = c/(1 - |z|/r) < \infty$.  Therefore $\| \bfU(\omega, z) \| < \infty$ for almost all $\omega \in \Omega$.  See \citet[Lemma 7.1, pp 182\textendash 183]{bos1}.  Clearly $\bfU(z) = {\mathcal M}[ \bfu ](z) = {\mathcal M}[ \bfu_+](z) = \bfU_+(z)$ because the transform uses no information for $t < 0$. If the time series is only observed for $t \in {\mathbb N}-1$ the ${\mathcal M}$-transform captures all of the observed information. Because $\bfu_+(t - s) = 0$ for $s > t$ it follows from (\ref{e:(A.1)}) and (\ref{e:(A.2)}) in Section~\ref{s:nt} of the Appendix that
\begin{equation}
\label{e:(3.9)}
\Delta^{-k} \bfu_+(t) = \mbox{$\sum_{s \in {\mathbb N}-1}$} \mbox{$\binom{k+s-1}{s}$} \bfu_+(t-s) = \mbox{$\sum_{s = 0}^t$} \mbox{$\binom{k+s-1}{s}$} \bfu(t-s)
\end{equation}
for $k \in {\mathbb N}$ and
\begin{equation}
\label{e:(3.10)}
\Delta^{\ell} \bfu_+(t) = \mbox{$\sum_{s=0}^{\ell}$} \mbox{$\binom{\ell}{s}$} (-1)^s \bfu_+(t-s) = \mbox{$\sum_{s=0}^{\min \{\ell,t\}}$} \mbox{$\binom{\ell}{s}$} (-1)^s \bfu(t-s)
\end{equation} 
for $\ell \in {\mathbb N}-1$ and all $t \in {\mathbb N}-1$. Therefore
\begin{eqnarray}
\label{e:(3.11)}
\lefteqn{{\mathcal M}[ \{ \Delta^{-k} \bfu_+(t) \}_{t \in {\mathbb N}-1}](z) = \mbox{$\sum_{t=0}^{\infty}$} \left[\, \mbox{$\sum_{s = 0}^t$} \mbox{$\binom{k+s-1}{s}$} \bfu(t-s) \right] z^t } \hspace{5cm} \nonumber \\
& & = \mbox{$\sum_{s=0}^{\infty}$} \mbox{$\binom{k+s-1}{s}$} z^s \left[\, \mbox{$\sum_{\tau=0}^{\infty}$} \bfu(\tau) z^{\tau} \right] = (1 - z)^{-k} \bfU(z) \hspace{1cm}
\end{eqnarray}
for $z \in D_r(0)$ and each $k \in {\mathbb N}$ and
\begin{eqnarray}
\label{e:(3.12)}
\lefteqn{{\mathcal M}[ \{ \Delta^{\ell} \bfu_+(t) \}_{t \in {\mathbb N}-1}](z) = \mbox{$\sum_{t=0}^{\infty}$} \left[\, \mbox{$ \sum_{s = 0}^{\min\{\ell,t \}}$} \mbox{$\binom{\,\ell \,}{s}$} (-1)^s \bfu(t-s) \right] z^t } \hspace{5cm} \nonumber \\
& = & \mbox{$\sum_{s=0}^{\ell}$} \mbox{$\binom{\, \ell \,}{s}$} (-1)^s z^s \left[\, \mbox{$\sum_{\tau=0}^{\infty}$} \bfu(\tau) z^{\tau} \right] = (1 - z)^{\ell} \bfU(z) \hspace{1cm}
\end{eqnarray}
for $z \in D_r(0)$ and each $\ell \in {\mathbb N}-1$.  The ${\mathcal M}$-transform satisfies a standard convolution identity. Let $\{V_t\}_{t \in {\mathbb N}-1} \in {\mathcal B}(Y,X)^{{\mathbb N}-1}$ with $\|V_t\| \leq c/r^t$ for some $c, r > 0$.  The corresponding ${\mathcal M}$-transform $S(z) = \sum_{t \in {\mathbb N}-1} V_t z^t \in {\mathcal B}(Y,X)$ is well defined for all $z \in D_r(0)$. Suppose $\bxi \in Y^{\mathbb Z}$ is a strong white noise with ${\mathcal M}$-transform $\bXi(z) = \sum_{j \in {\mathbb N}-1} \bxi(t)x^t$ for $z \in D_1(0)$.  If we define the convolution $\bfw = V \star \bxi \Leftrightarrow \bfw(t) = V_t \star \bxi(t)$ by setting
\begin{equation}
\label{e:(3.13)}
\bfw(t) = \mbox{$\sum_{s=0}^t$}\, V_s\, \bxi(t-s)
\end{equation}
for all $t \in {\mathbb N}-1$ then $\bfW(z) = S(z) \bXi(z) \Leftrightarrow {\mathcal M}[V \star \bxi](z) = {\mathcal M}[V](z) {\mathcal M}[\bxi](z)$ for all $z \in D_r(0) \cap D_1(0)$.
 
\section{Resolution of the \boldmath $AR(1)$ time series model}
\label{s:rtsm}

We consider the $AR(1)$ process defined in (\ref{e:(2.1)}) by $A_0 \bfx(t) + A_1 \bfx(t-1) = \bxi(t)$ for each $t \in {\mathbb N}-1$ with $\bfx(-1) = \bfc$.  It follows from this infinite collection of equations that
\begin{eqnarray*}
A_0 \bfx(0) + A_1 \bfx(-1) & = & \bxi(0) \\
A_0 \bfx(1)z + A_1 \bfx(0)z & = & \bxi(1)z \\
A_0 \bfx(2)z^2 + A_1 \bfx(1) z^2 & = & \bxi(2) z^2 \\
\vdots & = & \vdots
\end{eqnarray*}
for each $z \in {\mathbb C}$.  By summing the equations we can collect all observable information into a single algebraic equation
\begin{equation}
\label{e:(4.1)}
A(z) \bfX(z) = \bXi(z) - B_1\, \bfc
\end{equation}
where $A(z) = A_0 + A_1z$ is the autoregressive polynomial, $\bfX(z) = {\mathcal M}[\bfx](z)$ and $\bXi(z) = {\mathcal M}[\bxi](z)$ are the respective ${\mathcal M}$-transforms of $\bfx$ and $\bxi$ and where we prefer to write $A_1 = B_1$ on the right-hand side.   According to Lemma~\ref{lem:rcon} we may assume that $R(z) = A(z)^{-1} \in {\mathcal B}(Y,X)$ is analytic on a set $D_{1+\epsilon}(0) \setminus \{1\}$ for some $\epsilon > 0$ and has an isolated singularity at $z = 1$. An intuitive resolution of (\ref{e:(4.1)}) is therefore given by
\begin{equation}
\label{e:(4.2)}
\bfX(z) = R(z)[\bXi(z) - B_1\bfc] = R_{\mbox{\scriptsize \rm sin}}(z) [\bXi(z) - B_1\bfc] + R_{\mbox{\scriptsize \rm reg}}(z)[\bXi(z) - B_1\bfc]
\end{equation}
for all $z \in D_{1 + \epsilon}(0) \setminus \{1\}$. Our main task is to justify this intuitive resolution and obtain a suitable representation for $\bfx(t) = {\mathcal M}^{-1}[ \bfX](t)$ for $t \in {\mathbb N}-1$.

\subsection{A transform formula for the singular part of the resolvent}
\label{ss:tfspr}

The function $R_{\mbox{\scriptsize \rm sin}}(z)$ defined in (\ref{e:(3.6)}) is analytic for all $z \neq 1$. Therefore we can write
\begin{equation}
\label{e:(4.3)}
R_{\mbox{\scriptsize \rm sin}}(z) = \mbox{$\sum_{s \in {\mathbb N}-1}$} U_tz^t = {\mathcal M}[\{U_t\}_{t \in {\mathbb N}-1}](z)
\end{equation}
as an ${\mathcal M}$-transform for all $z \in D_{\, 1}(0)$ where the coefficients $\{U_t\}_{t \in {\mathbb N}-1}$ are given by the formula
\begin{equation}
\label{e:(4.4)}
U_t = (1/t!) \left[ d^{\, t} R_{\mbox{\scriptsize \rm sin}}(z)/dz^t \right]_{z=0} = (-1)(I_X - T_{-1}B_0)^{-t-1}T_{-1}.
\end{equation} 
Consequently $U_t \in {\mathcal B}(Q(Y), P(X))$ for all $t \in {\mathbb N}-1$. 

\subsection{A transform formula for the extended regular part of the resolvent}
\label{ss:tferpr}

The resolvent $R(z) = (A_0 + A_1z)^{-1}$ is analytic for $z \in D_{\, 1}(0)$. Therefore $A_0^{-1}$ is well defined and
\begin{equation}
\label{e:(4.5)}
R(z) = \mbox{$\sum_{t \in {\mathbb N}-1}$} R_t z^t = {\mathcal M}[\{R_t\}_{t \in {\mathbb N}-1}](z)
\end{equation}
for all $z \in D_{\, 1}(0)$ where
\begin{equation}
\label{e:(4.6)}
R_t = (1/t!) \left[ d^{\, t} (A_0 + A_1z)^{-1}/dz^t \right]_{z=0} = (-1)^t(A_0^{-1}A_1)^t A_0^{-1}
\end{equation}
for all $t \in {\mathbb N}-1$. It follows that $\lim_{t \rightarrow \infty} \|R_t \|^{1/t} = 1$.  We can extend the function $R_{\mbox{\scriptsize \rm reg}}(z)$ defined in (\ref{e:(3.7)}) to an analytic function $W(z)$ on $D_{1 + \epsilon}(0)$ by defining
\begin{equation}
\label{e:(4.7)}
W(z) = \left\{ \begin{array}{ll}
R(z) - R_{\mbox{\scriptsize \rm sin}}(z) & \mbox{for}\ z \in D_{1 + \epsilon}(0) \setminus \{1\} \\
T_0 & \mbox{for}\ z=1. \end{array} \right.
\end{equation}
We can now express the extended function $R_{\mbox{\scriptsize \rm reg}}(z) = W(z)$ as an ${\mathcal M}$-transform by writing
\begin{equation}
\label{e:(4.8)}
W(z) = \mbox{$\sum_{t \in {\mathbb N}-1}$} W_t z^t = {\mathcal M}[\{W_t\}_{t \in {\mathbb N}-1}](z)
\end{equation}
for $z \in D_{1+\epsilon}(0)$ where $W_t \in {\mathcal B}(Y,X)$ and $\lim_{t \rightarrow \infty} \|W_t\|^{1/t} \leq 1/(1 + \epsilon)$. Now (\ref{e:(4.7)}) gives
\begin{equation}
\label{e:(4.9)}
W_t = R_t - U_t = (-1)^t (A_0^{-1}A_1)^t A_0^{-1} + (I_X - T_{-1}B_0)^{-t-1}T_{-1}
\end{equation}
for each $t \in {\mathbb N}-1$. The inequality $\lim_{t \rightarrow \infty} \|W_t\|^{1/t} \leq 1/(1+ \epsilon)$ means that $\|W_t\| \leq c /(1 + \epsilon)^t$ for some $c > 0$ and all $t \in {\mathbb N}-1$. If we write $\zeta = z-1$ and recall that $R_{\mbox{\scriptsize \rm reg}}(z) = W(z)$ then
\begin{eqnarray}
\label{e:(4.10)}
\lefteqn{ R_{\mbox{\scriptsize \rm reg}}(z) = \mbox{$\sum_{t \in {\mathbb N}-1}$} W_t (1 + \zeta)^t = \mbox{$\sum_{t \in {\mathbb N}-1}$} W_t \left[\, \mbox{$\sum_{\ell=0}^t$} \mbox{$\binom{t}{\ell}$} \zeta^{\ell}\, \right] } \hspace{2cm} \nonumber \\
& & = \mbox{$\sum_{\ell \in {\mathbb N}-1}$} \left[\, \mbox{$\sum_{t \in {\mathbb N} - 1 + \ell}$} \mbox{$\binom{t}{\ell}$} W_t\, \right] \zeta^{\ell} = \mbox{$\sum_{\ell \in {\mathbb N}-1}$} \left[\, \mbox{$\sum_{r \in {\mathbb N}-1}$} \mbox{$\binom{\ell+r}{\ell}$} W_{\ell + r}\, \right] \zeta^{\ell} \hspace{1cm}
\end{eqnarray}
for all $z \in D_{\, \epsilon}(1)$. It follows from (\ref{e:(4.10)}) that we can also calculate $\{T_{\ell}\}_{\ell \in {\mathbb N}-1}$ using the formula
\begin{equation}
\label{e:(4.11)}
T_{\ell} = \mbox{$\sum_{r \in {\mathbb N}-1}$} \mbox{$\binom{\ell+r}{\ell}$} W_{\ell+r}
\end{equation}
for each $\ell \in {\mathbb N}-1$. To solve Example~\ref{ex:1} in subsection~\ref{ss:tferpr} we define $W(z) = R(z) - T_{-1}/(z-1)$ for $z \in D_{1 + \epsilon}(0)$ and show that $W_t \in {\mathcal B}(Q^c(Y),P^c(X))$ for all $t \in {\mathbb N}-1$. Lemma~\ref{lem:mtrr} shows that this is true in general. The proof is given in Section~\ref{s:pmr} of the Appendix.

\begin{lem}
\label{lem:mtrr}
Let $R_{\mbox{\scriptsize \rm reg}}(z) = W(z)$ be represented by the ${\mathcal M}$-transform $W(z) = \sum_{t \in {\mathbb N}-1} W_t z^t$ for all $z \in D_{1 + \epsilon}(0)$. Then $W_t = P^cW_t = W_tQ^c \in{\mathcal B}(Q^c(Y), P^c(X))$ for all $t \in {\mathbb N}-1$.
\end{lem}

We have shown that $R_{\mbox{\scriptsize \rm sin}}(z) = \sum_{t \in {\mathbb N}-1} U_t z^t = {\mathcal M}[\{U_t\}_{t \in {\mathbb N}-1}](z)$ and $R_{\mbox{\scriptsize \rm reg}}(z) = \sum_{t \in {\mathbb N}-1} W_t z^t = {\mathcal M}[\{W_t\}_{t \in {\mathbb N}-1}](z)$. We can also see that $\bXi(z) - B_1 \bfc = {\mathcal M}[\bxi(t) - B_1\bfc\, \bdelta(t)](z)$ where the series $\{ \bdelta(t) \}_{t \in {\mathbb Z}}$ is defined by $\bdelta(0)=1$ and $\bdelta(t) = 0$ if $t \neq 0$. It follows from (\ref{e:(3.13)}) and (\ref{e:(4.2)}) that
\begin{equation}
\label{e:(4.12)}
\bfx(t) = \bfx_{\mbox{\scriptsize \rm sin}}(t) + \bfx_{\mbox{\scriptsize \rm reg}}(t) = U_t \star [\bxi(t) - B_1 \bfc\, \bdelta(t)] + W_t \star [\bxi(t) - B_1\bfc\, \bdelta(t)]
\end{equation}
for all $t \in {\mathbb N}-1$. Formula (\ref{e:(4.12)}) is an implicit version of the extended GJRT.  It remains to show that (\ref{e:(4.12)}) can be rewritten explicitly in the form (\ref{e:(4.13)}) or in the modified form (\ref{e:(4.17)}).  The relevant proofs are given in Section~\ref{s:pmr} of the Appendix.

\begin{prop}[GJRT]
\label{p:1}
Under Assumption \ref{a:1} the series $\bfx(t)$ in \eqref{e:(2.1)} admits the solution
\begin{eqnarray}
\label{e:(4.13)}
\bfx(t) & = & \bfx_{\mbox{\scriptsize \rm sin}}(t) + \bfx_{\mbox{\scriptsize \rm reg}}(t) \nonumber \\
& = & \left[ \mbox{$\sum_{k \in {\mathbb N}}$} (-1)^k T_{-k} \Delta^{-k} \bxi_+(t) - U_t B_1 \bfc \right] +\, \left[ \rule{0cm}{0.45cm} \mbox{$\sum_{s=0}^t$} W_s\, \bxi(t-s) - W_tB_1 \bfc \right],
\end{eqnarray}
for all $t \in {\mathbb N}-1$, where $\bfx(-1) = \bfc$, $\bfx_{\mbox{\scriptsize \rm sin}}(t) = P\bfx(t) \in P(X)$, $\bfx_{\mbox{\scriptsize \rm reg}}(t) = P^c\bfx(t) \in P^c(X)$, and
$$
\begin{array}{l}
T_{-k}=(-1)^{k-1} (T_{-1}B_0)^{k-1}T_{-1}, \\
U_t=(-1)(I_X - T_{-1}B_0)^{-t-1}T_{-1},\\
W_t=(-1)^t (A_0^{-1}A_1)^t A_0^{-1} - U_t.
\end{array}
$$
\end{prop}

The following consequence of the GJRT is proved in Section~\ref{s:pmr} of the Appendix.

\begin{coro}
\label{c:1}
Under Assumption \ref{a:1} the series $\bfx \in X^{\mathbb Z}$ satisfies (\ref{e:(2.1)}) with $\bfx(-1) = \bfc$ if and only if $\bfx(t) = \bfx_{\mbox{\scriptsize \rm sin}}(t) + \bfx_{\mbox{\scriptsize \rm reg}}(t)$ where the singular part $\bfx_{\mbox{\scriptsize \rm sin}} \in P(X)^{{\mathbb N}-1}$ satisfies
\begin{equation}
\label{e:(4.14)}
QA_0P\bfx_{\mbox{\scriptsize \rm sin}}(t) + QA_1P \bfx_{\mbox{\scriptsize \rm sin}}(t-1) = Q\bxi(t)
\end{equation}
for $t \in {\mathbb N}-1$ and $\bfx_{\mbox{\scriptsize \rm sin}}(-1) = P\bfc$ and the regular part $\bfx_{\mbox{\scriptsize \rm reg}} \in P^c(X)^{{\mathbb N}-1}$ satisfies
\begin{equation}
\label{e:(4.15)}
Q^cA_0P^c\bfx_{\mbox{\scriptsize \rm reg}}(t) + Q^cA_1P^c \bfx_{\mbox{\scriptsize \rm reg}}(t-1) = Q^c\bxi(t)
\end{equation}
for $t \in {\mathbb N}-1$ and $\bfx_{\mbox{\scriptsize \rm reg}}(-1) = P^c\bfc$.  Therefore an integrated $AR(1)$ process $\bfx$ is the sum of separate $AR(1)$ processes on complementary subspaces.
\end{coro}

We define the extended singular and regular parts $\bfx_{\mbox{\scriptsize \rm sin}} \in P(X)^{\mathbb Z}$ and $\bfx_{\mbox{\scriptsize \rm reg}} \in P^c(X)^{\mathbb Z}$ by supposing that (\ref{e:(4.14)}) and (\ref{e:(4.15)}) are true for all $t \in {\mathbb Z}$.  It is convenient to define $A_{i, {\rm r}} = Q^cA_iP^c \in {\mathcal B}(P^c(X),Q^c(Y))$ for each $i=0,1$.   The resolvent $R_{\mbox{\scriptsize \rm reg}}(z) = (A_{0,{\rm r}} + A_{1,{\rm r}}z)^{-1} \in {\mathcal B}(Q^c(Y),P^c(X))$ is well defined and analytic for $z \in D_{1 + \epsilon}(0)$.  Therefore 
\begin{enumerate}
\item $A_{0,{\rm r}}^{-1} \in {\mathcal B}(Q^c(Y),P^c(X))$ is well defined,
\item $\Phi = (-1) A_{0,{\rm r}}^{-1}A_{1,{\rm r}} \in {\mathcal B}(P^c(X))$ is well defined, and
\item $\Vert \Phi^s \Vert = \Vert (A_{0,{\rm r}}^{-1}A_{1,{\rm r}})^s \Vert < c/(1 + \epsilon)^s$ for some $c > 0$ and all $s \in {\mathbb N}$.
\end{enumerate}
By rearranging the extended form of (\ref{e:(4.15)}) it follows that $\bfx_{\mbox{\scriptsize \rm reg}} \in P^c(X)^{\mathbb Z}$ is an $AR(1)$ process defined by the equation
\begin{equation}
\label{e:(4.16)}
\bfx_{\mbox{\scriptsize \rm reg}}(t) = \Phi \bfx_{\mbox{\scriptsize \rm reg}}(t-1) + \bfeta(t)
\end{equation}
for all $t \in {\mathbb Z}$ where $\bfx_{\mbox{\scriptsize \rm reg}}(-1) = P^c \bfc$ and where the strong white noise $\bfeta \in P^c(X)^{\mathbb Z}$ is defined by setting $\bfeta(t) = A_{0,{\rm r}}^{-1} Q^c \bxi(t)$ for each $t \in {\mathbb Z}$.  Our next result, Corollary~\ref{c:2}, is an alternative version of the extended GJRT.  This corollary is proved in Section~\ref{s:pmr} of the Appendix.

\begin{coro}
\label{c:2}
Under Assumption~\ref{a:1} the $AR(1)$ process $\bfx \in X^{\mathbb Z}$ satisfies
\begin{eqnarray}
\label{e:(4.17)}
\bfx(t) & = & \bfx_{\mbox{\scriptsize \rm sin}}(t) + \bfx_{\mbox{\scriptsize \rm reg}}(t) \nonumber \\
& = & \left[ \mbox{$\sum_{k \in {\mathbb N}}$} (-1)^k T_{-k} \Delta^{-k} \bxi_+(t) - U_t B_1 \bfc \right] + \left[ \rule{0cm}{0.45cm}\Phi^{t+1} \bfz_{\infty} + \mbox{$\sum_{s \in {\mathbb N}-1}$} \Phi^s \bfeta(t-s) \right]
\end{eqnarray}
for $t \in {\mathbb N}-1$ where $\bfz_{\infty} = \lim_{r \rightarrow \infty} \bfz_r = \lim_{r \rightarrow \infty} \Phi^r \bfx_{\mbox{\scriptsize \rm reg}}(-1-r) = P^c \bfc - \sum_{s \in {\mathbb N}-1} \Phi^s \bfeta(-1-s) \in P^c(X)$ almost surely.  If $P \bfc = \bfzero$ then $\bfx_{\mbox{\scriptsize \rm sin}} \sim I(d)$ for some $d \in {\mathbb N} \cup \{ \infty\}$.  If $\bfz_{\infty} = \bfzero$ then $\bfx_{\mbox{\scriptsize \rm reg}} \sim I(0)$.  If both conditions hold then $\bfx = (\bfx_{\mbox{\scriptsize \rm sin}} + \bfx_{\mbox{\scriptsize \rm reg}}) \sim I(d)$.
\end{coro} 

If $P^c \bfc = \sum_{s \in {\mathbb N}-1} \Phi^s \bfeta(-1-s)$ then $\bfz_{\infty} = \bfzero$ and $\bfx_{\mbox{\scriptsize \rm reg}}(t) = \sum_{s \in {\mathbb N}-1} \Phi^s \bfeta(t-s)$ is strictly stationary.   In this case $\bfx_{\mbox{\scriptsize \rm reg}}(t)$ is a standard linear process.  If $\bfz_{\infty} \neq \bfzero$ then ${\mathbb E}[ \| \Phi^{t+1}\bfz_{\infty}\| ] \rightarrow \bfzero$ as $t \rightarrow \infty$ and we say that $\bfx_{\mbox{\scriptsize \rm reg}}(t)$ is asymptotically strictly stationary.

\begin{rem}
\label{r:1}
Suppose $P^c(X) \neq \{ \bfzero \} \Leftrightarrow P^c \neq \bigzero$.  Choose $\bfp \in P^c(X)$ with $\bfp \neq \bfzero$ and define a bounded linear functional $\bff_{\sbfp} \neq \bfzero$ on the subspace $\{ \bfx \mid \bfx = \alpha \bfp\ \mbox{for all}\ \alpha \in {\mathbb C}\} \subset P^c(X)$ by setting $\langle \alpha \bfp, \bff_{\sbfp} \rangle = \alpha$ for all $\alpha \in {\mathbb C}$. The Hahn\textendash Banach theorem shows that we can extend $\bff_{\sbfp}$ to a bounded linear functional $\bff_{\sbfp} \in [P^c(X)]^*$. Now we can further extend $\bff_{\sbfp}$ to a bounded linear functional $\bff_{\sbfp} \in X^*$ by defining $\langle \bfx, \bff_{\sbfp} \rangle = \langle P^c\bfx, \bff_{\sbfp} \rangle$ for all $\bfx \in X$. Therefore $\langle P\bfx, \bff_{\sbfp} \rangle = 0$ for all $\bfx \in X$. Hence $\bff_{\sbfp} \in P(X)^{\perp}$. By Theorem 1 in \cite{alb5} one has $T_{-1}(Y) = P(X)$. Therefore $\bff_{\sbfp} \in T_{-1}(Y)^{\perp}$.  If $\bfx(-1) = \bfc = P \bfd + \sum_{s \in {\mathbb N}-1} \Phi^s \bfeta(-1-s)$ for some $\bfd \in X$ then $P^c \bfc = \sum_{s \in {\mathbb N}-1} \Phi^s \bfeta(-1-s)$ and so $\bfz_{\infty} = \bfzero$.  Now (\ref{e:(4.17)}) shows that $\langle \bfx(t), \bff_{\sbfp} \rangle = \sum_{s \in {\mathbb N}-1} \Phi^s \langle \bfeta(t-s), \bff_{\sbfp} \rangle$ for all $t \in {\mathbb N}-1$. Therefore $\langle \bfx, \bff_{\sbfp} \rangle$ is a standard linear process and $\langle \bfx, \bff_{\sbfp} \rangle \sim I(0)$. We conclude that the condition $P^c \neq \bigzero$ is sufficient to ensure that $\bfx$ is cointegrated.
\end{rem}

\begin{rem}
\label{r:2}
The models (\ref{e:(4.14)}) and (\ref{e:(4.15)}) describe $AR(1)$ processes on complementary subspaces. In an algebraic sense these processes are linearly independent. Any stochastic dependence is entirely due to stochastic dependence in the noise terms. If $Q\bxi$ and $Q^c\bxi$ are stochastically independent then so too are $\bfx_{\mbox{\scriptsize \rm sin}}$ and $\bfx_{\mbox{\scriptsize \rm reg}}$.
\end{rem}

\subsection{The attractor and cointegrating spaces}
\label{ss:acs}

We suppose Assumption~\ref{a:1} is true. In our notation the usual definition of the attractor space is $P(X) = T_{-1}(Y) \subset X$ while the usual definition of the cointegrating space is $P(X)^{\perp} = T_{-1}(Y)^{\perp} \subset X^*$ where $X^*$ is the dual space to $X$. The rationale for the definition of the attractor space\footnote{This \textit{rationale} was incorrectly stated in the original article: Phil Howlett, Brendan K. Beare, Massimo Franchi, John Boland and Konstantin Avrachenkov (2025), The Granger-Johansen representation theorem for integrated time series on Banach space,  \textit{Journal of Time Series Analysis}, \textbf{46}, 432\textendash 457, http://doi.org/ 10.1111/jtsa.12766.} is that each $P^c \bfx \in P^c(X)$ defined by $P^c \bfx(t) = \Phi^{t+1} \bfz_{\infty} + \sum_{s \in {\mathbb N}-1} \Phi^s \bfeta(t-s)$ for some $\bfz_{\infty} \in P^c(X)$ is a stable process with ${\mathbb E}[\| P^c \bfx(t) - \bfu_0(t)\|] \leq \| \Phi \|^{t+1} {\mathbb E}[ \|\bfz_{\infty} \|] \rightarrow 0$ as $t \rightarrow \infty$ where the time series $\bfu_0 \in P^c(X)$ is defined by $\bfu_0(t) = \sum_{s \in {\mathbb N}-1} \Phi^s \eta(t-s)$.  This means that deviations from the space $P(X)$ are restricted in the long run to deviations caused by the \textit{strictly stationary} standard linear process $\bfu_0 \sim I(0)$.  On the other hand ${\mathbb E}[ \| P \bfx \|] \rightarrow \infty$ as $t \rightarrow \infty$.  Thus the expected deviation relative to the expected norm converges to zero.  The rationale for the definition of the cointegrating space is that if $P^c(X) \neq \{ \bfzero\}$ then our earlier remarks show that there is some non-zero bounded linear functional $\bff_0 \in T_{-1}(Y)^{\perp} \subset X^*$ with $\langle \bfx, \bff_0 \rangle \sim I(0)$.

If $R(z)$ has a pole of order $d \in {\mathbb N}$ at $z=1$ there is a whole hierarchy of spaces with weaker versions of the above properties. If $\bfx \in T_{-k}(Y)$ then we can find $\bfy \in Y$ such that $\bfx = T_{-k} \bfy = T_{-k+1} [(-1) B_0T_{-1} \bfy] = T_{-k+1}\bfz$ where $\bfz = B_0T_{-1} \bfy \in Y$. Therefore $\bfx \in T_{-k+1}(Y)$. It follows that $T_{-1}(Y) \supset T_{-2}(Y) \supset \cdots \supset T_{-d}(Y) \supset T_{-d-1}(Y) = \{ \bfzero\}$. The inclusions are proper as the following argument shows. Suppose $T_{-r-1}(Y) = T_{-r}(Y) \neq \{ \bfzero \}$ for some $r \in {\mathbb N}$. Then $(-1)(T_{-1}B_0T_{-r-1}(Y) = (-1)(T_{-1}B_0T_{-r}(Y) \Leftrightarrow T_{-r-2}(Y) = T_{-r-1}(Y)$.  Therefore $T_{-r-2}(Y) = T_{-r}(Y)$.  An inductive argument now shows that $T_{-k}(Y) = T_{-r}(Y)$ for all $k \in {\mathbb N}+r$ but this is a contradiction because $T_{-d-1} = \{ \bfzero\}$.  A straightforward argument also shows that $T_{-1}(Y)^{\perp} \subset T_{-2}(Y)^{\perp} \subset \cdots \subset T_{-d}(Y)^{\perp} \subset X^*$.  Once again the inclusions are proper.  If $\bff \in T_{-r}(Y)^{\perp} \setminus T_{-r +1}(Y)^{\perp}$ for some $r \in \{2,3,\ldots,d-1\}$ then $\bff \neq \bfzero$ and $\langle \bfx, \bff \rangle \sim I(r)$.

If $R(z)$ has an essential singularity at $z = 1$ the above inclusions are still valid but it is no longer necessary that the inclusions are proper. Nevertheless it is still true that if $\bff \in T_{-r}(Y)^{\perp} \setminus T_{-r +1}(Y)^{\perp}$ for some $r \in \{2,3, \cdots \}$ then $\bff \neq \bfzero$ and $\langle \bfx, \bff \rangle \sim I(r)$.

\section{Special cases: \boldmath The GJRT when $R(z)$ has a pole at $z=1$}
\label{s:sc}

\cite{bea2} consider specific results for $I(d)$ processes when $d=1,2$. We will show that their results are special cases of our extended results. We reiterate that our representation is also new for cointegrated $AR(1)$ vector processes on ${\mathbb R}^n$. We suppose only that $R(z)$ is analytic on a set $[D_{1}(0)^a \setminus \{1\}] \cup D_{0, \delta}(1)$ for some $\delta > 0$ and has an isolated singularity at $z = 1$.

\subsection{The GJRT when \boldmath $R(z)$ has a simple pole at $z=1$}
\label{ss:scd1}

If $R(z)$ has a simple pole at $z=1$ then $T_{-k} =(-1)^{k-1}(T_{-1}B_0)^{k-1}T_{-1} = \bigzero$ for $k > 1$ and the binomial expansion shows that
$$
U_t = (-1)[I_X - T_{-1}B_0]^{-t-1}T_{-1} = (-1)[I_X + (t+1)(T_{-1}B_0) + \cdots ]T_{-1} = (-1) T_{-1}
$$
for all $t \in {\mathbb N}-1$.  Therefore (\ref{e:(4.17)}) becomes
\begin{equation}
\label{e:(5.1)}
\bfx(t) = (-1) T_{-1} \mbox{$\sum_{s=0}^t$}\, \bxi(t-s) + T_{-1}B_1\bfc + \Phi^{t+1} \bfz_{\infty} +  \mbox{$\sum_{s \in {\mathbb N}-1}$}\, \Phi^s \bfeta(t-s) 
\end{equation}
for all $t \in {\mathbb N}-1$ where we have used the identity $\Delta^{-1} \bxi_+(t) = \sum_{s=0}^t \bxi(t-s)$ for all $t \in {\mathbb N}-1$.  Therefore $\bfx \not \sim I(0)$.

\begin{rem}
\label{r:3}
If $\bfc = P\bfd +\sum_{s \in {\mathbb N}-1} \Phi^s \bfeta(-1-s)$ for some $\bfd \in X$ then $\bfz_{\infty} = \bfzero$ and (\ref{e:(5.1)}) implies
$$
\Delta \bfx(t) = (-1)T_{-1} \bxi(t) + \mbox{$\sum_{s \in {\mathbb N}-1}$} \Phi^s \Delta \bfeta(t-s)
$$
for all $t \in {\mathbb N}$.  Therefore $\Delta \bfx = \{ \Delta \bfx(t) \}_{t \in {\mathbb N}} \in X^{\mathbb N}$ is a strictly stationary standard linear process.  Thus $\Delta \bfx \sim I(0)$.  It follows that $\bfx \sim I(1)$.  Our results provide a natural extension of the results obtained by \cite{bea2}.
\end{rem}

\subsection{The GJRT when \boldmath $R(z)$ has a pole of order $2$ at $z=1$}
\label{ss:scd2}

If $R(z)$ has a pole of order two at $z = 1$ then $T_{-k} = (-1)^{k-1}(T_{-1}B_0)^{k-1}T_{-1} = \bigzero$ for $k > 2$ and the binomial expansion shows that
$$
U_t = (-1)(I_X - T_{-1}B_0)^{-t-1}T_{-1} = (-1)[I_X + (t+1)(T_{-1}B_0) + \cdots ]T_{-1} = - T_{-1} + (t+1)T_{-2}
$$
for all $t \in {\mathbb N}-1$.  The identities $\Delta^{-1} \bxi_+(t) = \sum_{s=0}^t\, \bxi(t-s)$ and $\Delta^{-2} \bxi_+(t) = \sum_{s=0}^t\, (s+1) \bxi(t-s)$ can then be used to show that (\ref{e:(4.17)}) becomes
\begin{eqnarray}
\label{e:(5.2)}
\bfx(t) \!\!\! & = &\!\!\!  \left[ T_{-2} \mbox{$\sum_{s=0}^t$} (s+1) \bxi(t-s) - T_{-1} \mbox{$\sum_{s=0}^t$} \bxi(t-s) - (t+1)T_{-2}B_1\bfc + T_{-1}B_1 \bfc \right] \hspace{1cm} \nonumber \\
& & \hspace{5cm} +\left[ \Phi^{t+1} \bfz_{\infty} + \mbox{$\sum_{s \in {\mathbb N}-1}$}\, \Phi^s \bfeta(t-s) \right]
\end{eqnarray}
for all $t \in {\mathbb N}-1$.

\begin{rem}
\label{r:4}
If $\bfc = P\bfd + \sum_{s \in {\mathbb N}-1} \Phi^s \bfeta(-1-s)$ for some $\bfd \in X$ then $\bfz_{\infty} = \bfzero$ and (\ref{e:(5.2)}) becomes
\begin{eqnarray}
\label{e:(5.3)}
\lefteqn{\bfx(t) = T_{-2} \mbox{$\sum_{s=0}^t$} (s+1) \bxi(t-s) - T_{-1} \mbox{$\sum_{s=0}^t$} \bxi(t-s) } \hspace{4cm} \nonumber \\
& &  - (t+1)T_{-2}B_1\bfd + T_{-1}B_1 \bfd + \mbox{$\sum_{s \in {\mathbb N}-1}$}\, \Phi^s \bfeta(t-s)
\end{eqnarray}
for all $t \in {\mathbb N}-1$.  Therefore
$$
\Delta \bfx(t) = T_{-2} \mbox{$\sum_{s=0}^t$} \bxi(t-s) - T_{-1}\bxi(t) - T_{-2}B_1\bfd + \mbox{$\sum_{s \in {\mathbb N}-1}$}\, \Phi^s \Delta \bfeta(t-s)
$$
for all $t \in {\mathbb N}$ and
$$
\Delta^2 \bfx(t) = T_{-2} \bxi(t) - T_{-1}\Delta \bxi(t) + \mbox{$\sum_{s \in {\mathbb N}-1}$}\, \Phi^s \Delta^2 \bfeta(t-s)
$$
for all $t \in {\mathbb N}+1$.  Therefore $\Delta^2 \bfx = \{ \Delta^2 \bfx(t)\}_{t \in {\mathbb N}+1} \in X^{{\mathbb N}+1}$ is a strictly stationary standard linear process.  Therefore $\Delta^2 \bfx \sim I(0)$.  Clearly $\bfx \not \sim I(0)$ and $\Delta \bfx \not \sim I(0)$.  It follows that $\bfx \sim I(2)$.  Once again these results are a natural extension of the results reported in \cite{bea2}.
\end{rem}

\section{Three examples to illustrate the proposed representation}
\label{ss:tepr}

The first example is a two-dimensional time series where $R(z)$ has a simple pole at $z=1$.  We show that $\bfx$ is the sum of two separate $AR(1)$ time series.  The second example is an infinite-dimensional time series where $R(z)$ has an essential singularity at $z=1$.  The third example is taken from \citet[Examples 3.3, pp 757\textendash 758 and 3.4, pp 763\textendash 764]{seo1} and concerns an infinite-dimensional series where $R(z)$ has a pole of order two at $z=1$.   This example is used to demonstrate the general solution procedure on separable Banach space.  The general solution procedure is described in \cite{alb5}.

\begin{exa}[\boldmath $\bfx \sim I(1)$]
\label{ex:1}

{\rm Let $X = Y = {\mathbb R}^2$ and let $\epsilon \in (0,1)$. Let $\bfx \in X^{{\mathbb N}-1}$ be an $AR(1)$ process with $A_0 \bfx(t) + A_1 \bfx(t-1) = \bxi(t)$ for all $t \in {\mathbb N}-1$ where
$$
A_0 = \left[ \begin{array}{cc}
0 & -1 \\
1+\epsilon & -(2+\epsilon) \end{array} \right], \quad A_1 = \left[ \begin{array}{cc}
1 & 0 \\
0 & 1 \end{array} \right] \quad \mbox{and} \quad A(z) = \left[ \begin{array}{cc}
z & -1 \\
1+\epsilon & z-(2+\epsilon) \end{array} \right].
$$
We assume that $\bfx(-1) = \bfzero$ and that ${\mathbb E}[\bxi(t)] = \bfzero \in {\mathbb R}^2$ and ${\mathbb E}[ \| \bxi(t) \|^2] = \sigma^2$ for all $t \in {\mathbb N}-1$. We define $B_0 = A_0 + A_1$ and $B_1 = A_1$. The resolvent is given by
$$
R(z) = (z-1)^{-1} \left[ \begin{array}{cc}
1 + \epsilon^{-1} & - \epsilon^{-1} \\
1 + \epsilon^{-1} & - \epsilon^{-1} \end{array} \right] + (z-1-\epsilon)^{-1} \left[ \begin{array}{cc}
- \epsilon^{-1} & \epsilon^{-1} \\
-(1 + \epsilon^{-1}) & 1 + \epsilon^{-1} \end{array} \right]
$$
for all $z \in {\mathbb C} \setminus \sigma$ where $\sigma = \{1, 1+\epsilon \}$. Therefore $R(z)$ has a pole of order one at $z = 1$ and we can write $R(z) = T_{-1}(z-1)^{-1} + \sum_{\ell \in {{\mathbb N}-1}} T_{\ell}(z-1)^{\ell}$ as a Laurent series for $z \in D_{0,\epsilon}(1)$ where
$$
T_{-1} = \left[ \begin{array}{cc}
1 + \epsilon^{-1} & -\epsilon^{-1} \\
1 + \epsilon^{-1} & -\epsilon^{-1} \end{array} \right], \quad T_0 = \epsilon^{-1} \left[ \begin{array}{cc}
\epsilon^{-1} & -\epsilon^{-1} \\
(1 + \epsilon^{-1}) & -(1 + \epsilon^{-1}) \end{array} \right] \quad \mbox{and} \quad T_{\ell} = \epsilon^{-\ell} T_0.
$$
The spectral projection on $X$ is $P = T_{-1}B_1 = T_{-1}$ with $P(X) = \{ \bfx \mid P\bfx = \bfx \} = \mbox{sp}\{\bfe_1 + \bfe_2\}$ and the complementary projection $P^c = I_X - P$ with $P^c(X) = \{ \bfx \mid P^c\bfx = \bfx \} = \mbox{sp}\{\bfe_1 + (1 + \epsilon)\bfe_2\}$. The spectral projection on $Y$ is $Q = B_1T_{-1} = T_{-1}$ where similar remarks apply. The singular part of the resolvent $R_{\mbox {\scriptsize \rm sin}}(z) \in {\mathcal B}(Q(Y),P(X))$ can be written directly as an ${\mathcal M}$-transform
$$
R_{\mbox {\scriptsize \rm sin}}(z) = T_{-1}(z-1)^{-1} = \mbox{$\sum_{t \in {\mathbb N}-1}$} (-1)T_{-1} z^t = \mbox{$\sum_{t \in {{\mathbb N}-1}}$} U_tz^t 
$$
for $z \in D_1(0)$ where $U_t = (-1)T_{-1}$. The same direct approach for the regular part $R_{\mbox {\scriptsize \rm reg}}(z)$ gives
\begin{eqnarray*}
R_{\mbox {\scriptsize \rm reg}}(z) & = & \mbox{$\sum_{\ell \in {{\mathbb N}-1}}$} \epsilon^{-\ell} T_0 (z-1)^{\ell} = \mbox{$\sum_{\ell \in {{\mathbb N}-1}}$} (-1)^{\ell} \epsilon^{-\ell} T_0\, \mbox{$\sum_{t=0}^{\ell} \binom{\ell}{t}$} z^t \\
& &\hspace{6cm} =\, \mbox{$\sum_{t \in {{\mathbb N}-1}}$} \left[ \mbox{$\sum_{\ell \in {\mathbb N}+t-1} \binom{\ell}{t}$} (-1)^{\ell} \epsilon^{-\ell} \right] T_0\, z^t
\end{eqnarray*}
which is not valid because $\epsilon \in (0,1)$ and each coefficient of $z^t$ involves the sum of a divergent series. However $R_{\mbox{\scriptsize \rm reg}}(z) \in {\mathcal B}(Q^c(Y),P^c(X))$ can be written as an ${\mathcal M}$-transform using the extended definition
$$
R_{\mbox{\scriptsize \rm reg}}(z) = R(z) - R_{\mbox{\scriptsize \rm sin}}(z) = \left[1 - (1+ \epsilon)^{-1}z \right]^{-1} \epsilon (1+ \epsilon)^{-1} T_0 = \mbox{$\sum_{t \in {{\mathbb N}-1}}$} W_t\, z^t 
$$
for $z \in D_{1 +\epsilon}(0) \setminus \{1\}$ and $R_{\mbox{\scriptsize \rm reg}}(1) = T_0$ where $W_t = \epsilon (1 + \epsilon)^{-(t+1)} T_0 \in {\mathcal B}(Q^c(Y),P^c(X))$. We can use the relationship $R(z) = R_{\mbox{\scriptsize \rm sin}}(z) + R_{\mbox{\scriptsize \rm reg}}(z)$ and the initial condition $\bfx(-1) = \bfzero$ to rewrite (\ref{e:(4.2)}) as two equations $\bfX_{\mbox{\scriptsize \rm sin}}(z) = R_{\mbox{\scriptsize \rm sin}}(z) Q \bXi(z)$ and $\bfX_{\mbox{\scriptsize \rm reg}}(z) = R_{\mbox{\scriptsize \rm reg}}(z) Q^c \bXi(z)$ where $\bfX(z) = \bfX_{\mbox{\scriptsize \rm sin}}(z) + \bfX_{\mbox{\scriptsize \rm reg}}(z)$. Therefore $\bfx(t) = \bfx_{\mbox{\scriptsize \rm sin}}(t) + \bfx_{\mbox{\scriptsize \rm reg}}(t) = U_t \star [Q\bxi](t) + W_t \star [Q^c \bxi](t)$ for all $t \in {\mathbb N}-1$ where $\bfx_{\mbox{\scriptsize \rm sin}} \sim I(1)$ is an $AR(1)$ process in $P(X)$ satisfying
$$
(QA_0P) \bfx_{\mbox{\scriptsize \rm sin}}(t) + (QA_1P) \bfx_{\mbox{\scriptsize \rm sin}}(t-1) = Q \bxi(t)
$$
with $\bfx_{\mbox{\scriptsize \rm sin}}(-1) = \bfzero$ and $\bfx_{\mbox{\scriptsize \rm reg}} \sim I(0)$ is an $AR(1)$ process in $P^c(X)$ satisfying
$$
(Q^cA_0P^c) \bfx_{\mbox{\scriptsize \rm reg}}(t) + (Q^cA_1P^c) \bfx_{\mbox{\scriptsize \rm reg}}(t-1) = Q^c \bxi(t)
$$
with $\bfx_{\mbox{\scriptsize \rm reg}}(-1) = \bfzero$. The original time series $\bfx$ is therefore the sum of separate $AR(1)$ time series on complementary subspaces.  We show in Corollary~\ref{c:1} that this is true in general.} $\hfill \Box$
\end{exa}

\begin{figure}[htb]
\begin{center}
\includegraphics[width=12cm]{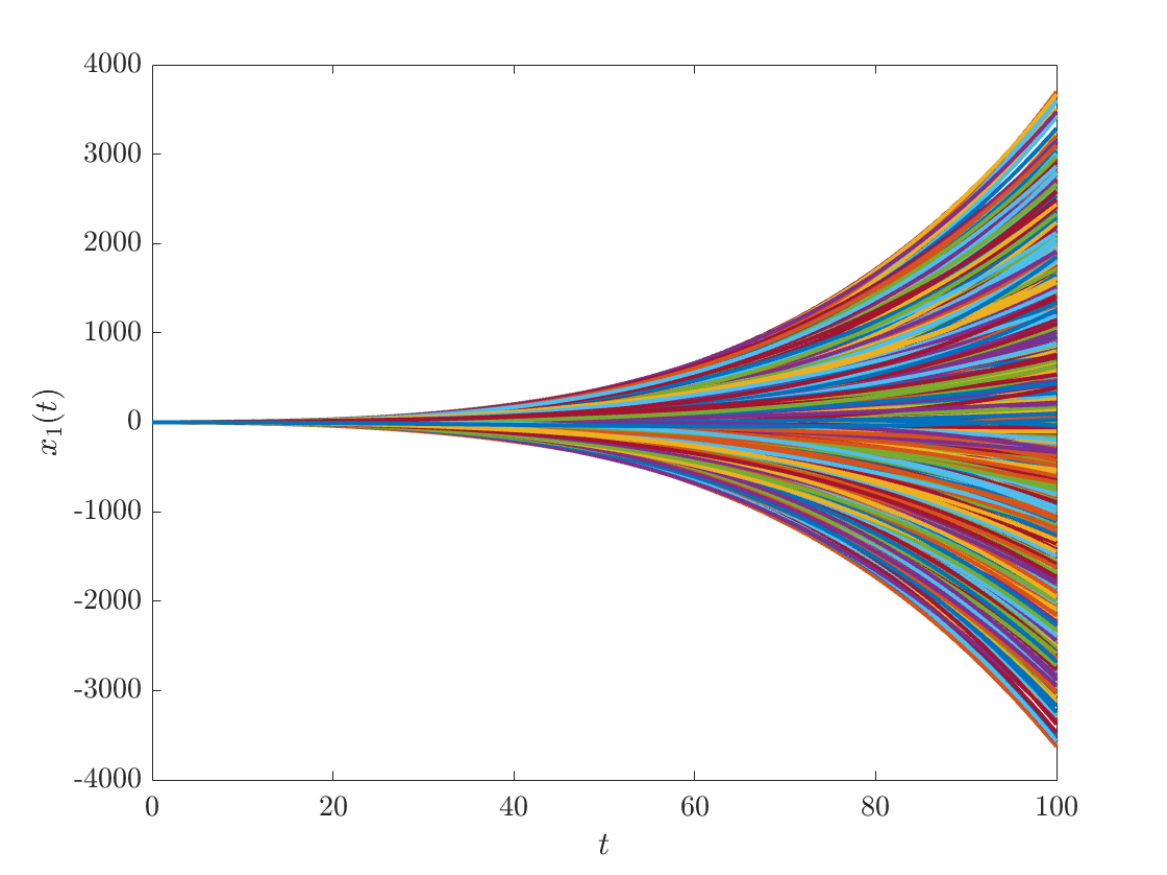}
\end{center}
\caption{\small Composite graph for the simulated values $\{\bfx_1(t)\}_{t \in \{0,1,\ldots,100\}}$ in Example~\ref{ex:2}.}
\label{fig1}
\end{figure}

\begin{figure}[htb]
\begin{center}
\includegraphics[width=12cm]{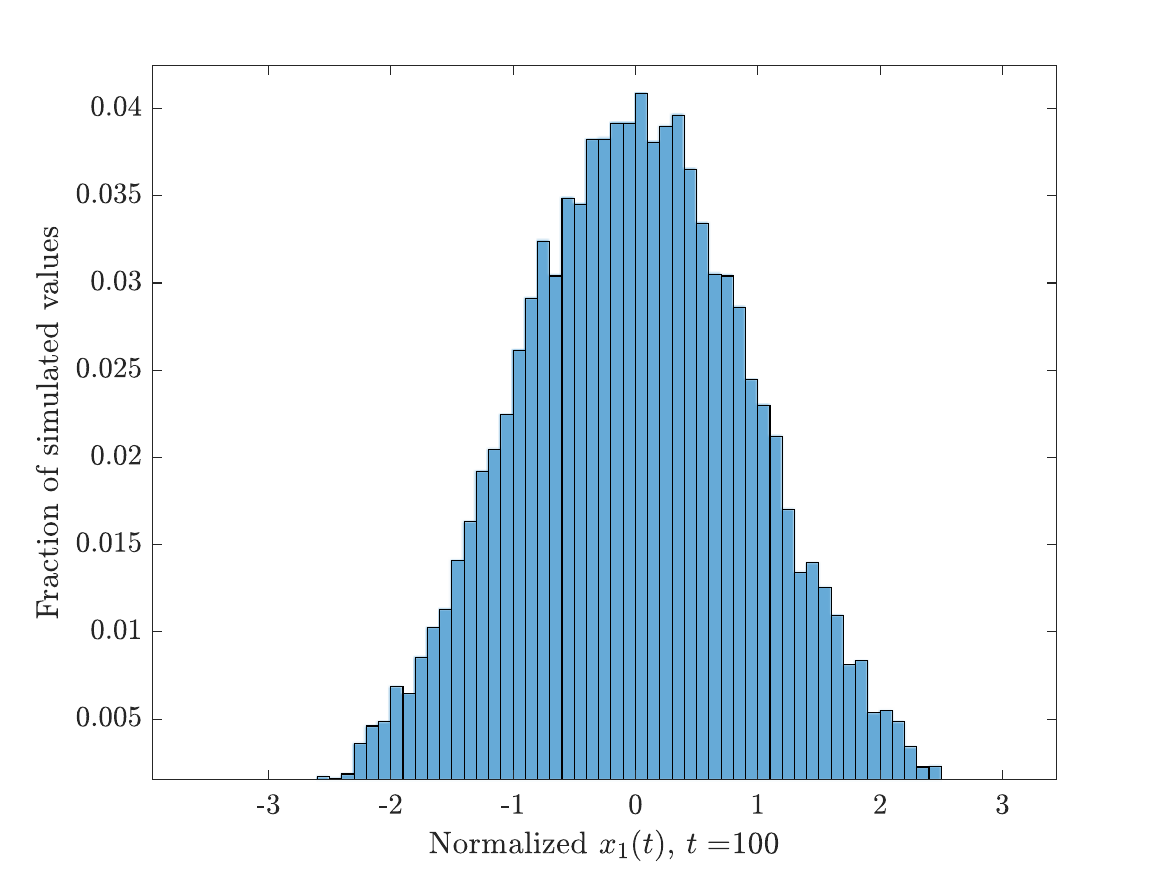}
\end{center}
\caption{\small Histogram showing the frequency distribution of simulated values for the normalized variable $\{\bfx_1(t) - m[\bfx_1(t)]\}/s[\bfx_1(t)]$ at $t =100$ in Example~\ref{ex:2}.}
\label{fig2}
\end{figure}

\begin{exa}[\boldmath $\bfx \sim I(\infty)$]
\label{ex:2}

{\rm Let $\{\lambda_j\}_{j \in {\mathbb N}} \subset {\mathbb R}$ with $\lambda_j > 0$ and $\lambda_j^{1/j} \downarrow 0$ as $j \uparrow \infty$ and $\{ \mu_j \}_{j \in {\mathbb N}} \subset {\mathbb R}$ with $\mu_j > 0$ and $\sum_{j \in {\mathbb N}} \mu_j = 1$. Suppose that $\bfx_j = \{ \bfx_j(t) \}_{t \in {\mathbb Z}} \in {\mathbb R}^{\mathbb Z}$ for each $j \in {\mathbb N}$ and that $\bfx_j, \bfx_{j+1}$ satisfy the recurrence relation
$$
\bfx_j(t) - \bfx_j(t-1) - \lambda_j \bfx_{j+1}(t) = \mu_j \bxi(t)
$$
with $\bfx_j(-1) = 0$ for each $j \in {\mathbb N}$ where $\bxi \in {\mathbb R}^{\mathbb Z}$ is a strong white noise. An ${\mathcal M}$-transform gives
$$
\bfX_j(z) - z\bfX_j(z) - \lambda_j \bfX_{j+1}(z) = \mu_j \bXi(z)
$$
for each $j \in {\mathbb N}$, where $\bfX_j(z) = {\mathcal M}[x_j](z)$ and $\bXi(z) = {\mathcal M}[\bxi](z)$.  In vector form we can write 
$$
[ (1 - z)I - U] \bfX(z) = \bXi(z) \bmu
$$
for each $z \in {\mathbb C}$ where $\bfX(z) = \sum_{j \in {\mathbb N}} \bfX_j(z) \bfe_j \in {\mathbb C}^{\infty}$,  $\bmu = \sum_{j \in {\mathbb N}} \mu_j \bfe_j \in \ell^1$, $I = [\bfe_1, \bfe_2, \ldots] \in {\mathcal B}(\ell^1)$ is the identity operator and $U = [ \bfzero, \lambda_1 \bfe_1, \lambda_2 \bfe_2, \lambda_3 \bfe_3, \ldots ] \in {\mathcal B}(\ell^1)$ is a strictly upper triangular operator. We can see that
$$
U^2 = [ \bfzero, \bfzero, \lambda_1\lambda_2 \bfe_1, \lambda_2 \lambda_3 \bfe_2, \lambda_3 \lambda_4 \bfe_4,\ldots ], \ \ U^3 = [ \bfzero, \bfzero, \bfzero, \lambda_1 \lambda_2 \lambda_3 \bfe_1, \lambda_2 \lambda_3 \lambda_4 \bfe_2, \lambda_3 \lambda_4 \lambda_5 \bfe_3, \ldots], \ \ \ldots
$$
and so on. Thus $U^n \bfe_k = \bfzero$ for $k \leq n$ and $U^n \bfe_{n+k} = \lambda_k \lambda_{k+1} \cdots \lambda_{n+k-1} \bfe_k$ for $n, k \in {\mathbb N}$. It follows that $\|U^n\|^{1/n} < \lambda_1 \lambda_n^{1/n} \downarrow 0$ as $n \uparrow \infty$. A Neumann expansion now shows that
$$
R(z) = [(1 - z)I - U]^{-1} = I/(1-z) + U/(1-z)^2 + U^2/(1-z)^3 + \cdots 
$$
for all $z \neq 1$. Therefore $R(z)$ has an isolated essential singularity at $z = 1$. Hence $\bfx \sim I(\infty)$. If we define $V = {\mathcal M}^{-1}[R]$ then
$$
V(t) = \left[ I + U \Delta^{-1} + U^2 \Delta^{-2} + \cdots \right] [\bfone] (t) = I + \mbox{$\binom{t+1}{1}$}U + \mbox{$\binom{t+2}{2}$}U^2 + \cdots
$$
where $\bfone(t) = 1$ and $\binom{t+n}{n} = (t+n)(t+n-1) \cdots (t+1)/n!$ for all $t \in {{\mathbb N}-1}$ and $n \in {\mathbb N}$. Because $\bXi(z) \in {\mathbb C}$ is a scalar it follows that
$$
{\mathcal M}^{-1}[R(z)\bXi(z)\bmu](t) = {\mathcal M}^{-1}[\bXi(z)R(z) \bmu](t)= (\bxi \ast V)(t) \bmu
$$
where $(\bxi \ast V)(t) = \sum_{s=0}^t \bxi(t-s)V(s)$ for all $t \in {\mathbb N}-1$. Therefore
$$
\bfx(t) = \mbox{$\sum_{n \in {\mathbb N} - 1}$} \left[ \mbox{$\sum_{s=0}^t$}\, \bxi(t - s) \mbox{$\binom{s+n}{n}$} \right] U^n \bmu
$$
for all $t \in {\mathbb N}-1$. We can extract the component $\bfx_k(t) = \langle \bfe_k, \bfx(t) \rangle$ for each $k \in {\mathbb N}$. In this regard we note that $\langle \bfe_k, U^n \bmu \rangle = \lambda_k \cdots \lambda_{k+n-1} \mu_{k+n}$ for $n \in {\mathbb N}$ and so
\begin{eqnarray*}
\bfx_k(t) & = & \left[ \mbox{$\sum_{s=0}^t$}\, \bxi(t-s) \right] \mu_k + \mbox{$\sum_{n \in {\mathbb N}}$} \left[ \mbox{$\sum_{s=0}^t$}\, \bxi(t - s) \mbox{$\binom{s+n}{n}$} \right] \lambda_k \cdots \lambda_{k+n-1} \mu_{k+n} \\
& = &\Delta^{-1} \bxi_+(t) \mu_k + \mbox{$\sum_{n \in {\mathbb N}}$} \left[ \Delta^{-n-1} \bxi_+(t) \right] \lambda_k \cdots \lambda_{k+n-1} \mu_{k+n}
\end{eqnarray*}
for all $t \in {\mathbb N}-1$. The final form shows that no finite order difference will eliminate all of the unbounded terms. This confirms that $\bfx_k \sim I(\infty)$ for each $k \in {\mathbb N}$.

We used M{\sc atlab} to simulate the evolution of the component $\bfx_1 \sim I(\infty)$ setting $\lambda_j = 1/j!$ and $\mu_j = 1/2^{j+1}$ with noise defined by a pseudo-random variable with a uniform distribution on the interval $(-1/2,1/2)$.  We performed $20,000$ independent trials for the series $\{ \bfx_1(t)\}_{t \in \{0,1,\ldots,100\}}$. The composite graph is shown in Figure~\ref{fig1}.  The histogram for the frequency distribution of the normalized values $\{\bfx_1(t) - m[\bfx_1(t)]\}/s[\bfx_1(t)]$ at $t =100$ is shown in Figure~\ref{fig2}. The mean and standard deviation for the trials were $m[x_1(100)] \approx -12.8811$ and $s[\bfx_1(100)] \approx 1.0051 \times 10^3$.  Although the composite graph shows that the range of values for $\bfx_1(t)$ increases rapidly as $t$ increases the histogram shows that the maximum frequency occurs when $\bfx_1(t) \approx 0$.} $\hfill \Box$
\end{exa}

\begin{exa}[\textbf{The general solution process}]
\label{ex:3}
{\rm Let $X = Y = c_0$ be the space of real-valued vectors $\bfx = [ \bfx_n ]_{n \in {\mathbb N}}$ with $\lim_{n \rightarrow \infty} \bfx_n = 0$ and norm $\| \bfx \|_{\infty} = \max_{n \in {\mathbb N}} |\bfx_n|$. The space $c_0$ is a closed separable subspace of $\ell^{\infty}$ with basis $\{\bfe_n\}_{n \in {\mathbb N}}$.  Let $\lambda \in (0,1)$ and suppose the time series $\bfx \in X^{\mathbb Z}$ satisfies $A_0 \bfx(t) + A_1\bfx(t-1) = \bxi(t)$ for all $t \in {\mathbb Z}$ with $\bfx(-1) = \bfc \in X$ where
$$
A_0 = I = [\bfe_1,\bfe_2,\ldots] \in {\mathcal B}(X,Y)  \quad \mbox{and} \quad A_1 = (-1)[\bfe_1, \bfe_1 + \bfe_2, \lambda \bfe_3,\lambda^2 \bfe_4, \lambda^3 \bfe_5,\ldots] \in {\mathcal B}(X,Y).
$$
The autoregressive polynomial is $A(z) = A_0 + A_1z \in {\mathcal B}(X,Y)$ for all $z \in {\mathbb C}$.  The operator $A(1) = A_0 + A_1$ is singular and so we write $A(z) = B_0 + B_1(z-1)$ where
$$
B_0 = A_0 + A_1 = [\bfzero, -\bfe_1, (1-\lambda)\bfe_3,(1- \lambda^2)\bfe_4, (1-\lambda^3)\bfe_5, \ldots] \quad \mbox{and} \quad B_1 = A_1.
$$
We wish to find a Laurent series $R(z) = \sum_{j \in {\mathbb Z}} T_j(z-1)^j$ which converges on some region $z \in D_{0,\epsilon}(1)$. Thus we need to find $\{T_j\}_{j \in {\mathbb Z}} \in {\mathcal B}(Y,X)$ satisfying the fundamental equations (\ref{e:(3.1)}) and (\ref{e:(3.2)}) and the magnitude constraints (\ref{e:(3.3)}). We follow \citet{alb5} and begin by finding the infinite-length Jordan chains at $z=1$.

The singular Jordan chains $\{ \bfu_{-n} \}_{n \in {\mathbb N}}$ satisfy
\begin{equation}
\label{e:(6.1)}
B_0 \bfu_{-n} + B_1 \bfu_{-n-1} = \bfzero
\end{equation}
for all $n \in {\mathbb N}$ with $\| \bfu_{-n} \|_{\infty}^{1/n} \rightarrow 0$ as $n \rightarrow \infty$. If we set $\bfu_{-1} = \bft = \bft_1\bfe_1 + \bft_2\bfe_2 + \bft_3 \bfe_3 + \cdots$ then (\ref{e:(6.1)}) shows that $\bfu_{-n,2} = 0$ for $n \geq 2$, $\bfu_{-n,1} = - \bfu_{-n+1,2} = 0$ for $n \geq 3$ and $\bfu_{-n, j} = (1 - 1/\lambda^{j-2})^{n-1} \bft_j$ for all $n \in {\mathbb N}$ and all $j \geq 3$. If $\bft_j \neq 0$ for $j \geq 3$ then
$$
\| \bfu_{-n} \|_{\infty}^{1/n} \geq (1/\lambda^{j-2} - 1)^{1 - 1/n} |\bft_j|^{1/n} \rightarrow (1/\lambda^{j-2} - 1) > 0
$$
for each $j \geq 3$ as $n \rightarrow \infty$. Thus we must choose $\bft_j = 0$ for $j \geq 3$. Therefore
$\bfu_{-1} = \bft_1 \bfe_1 + \bft_2 \bfe_2$, $\bfu_{-2} = \alpha \bfe_1$ for some $\alpha \in {\mathbb R}$ and $\bfu_{-n} = \bfzero$ for $n \geq 3$. This means the subspace $X_{\mbox{\scriptsize \rm sin}} \subset X$ is the closed subspace $\mbox{sp}(\{ \bfe_1, \bfe_2 \})\, \subset c_0$.

The regular Jordan chains $\{ \bfu_n\}_{n \in {\mathbb N}}$ satisfy
\begin{equation}
\label{e:(6.2)}
B_0 \bfu_{n+1} + B_1 \bfu_n = \bfzero
\end{equation}
for all $n \in {\mathbb N}$ with $\| \bfu_n \|_{\infty}^{1/n} \rightarrow a$ for some $a > 0$. If we set $\bfu_1 = \bfs = \bfs_1\bfe_1 + \bfs_2\bfe_2 + \bfs_3\bfe_3 + \cdots$ then (\ref{e:(6.2)}) shows that $\bfu_{1,2} = \bfs_2 = 0$. We also have $\bfu_{2,2} = 0$ and $\bfu_{2,2} = - \bfu_{1,1} = - \bfs_1$. Therefore $\bfs_1 = 0$ as well. The remaining equations show that $\bfu_{n,j} = [1  - 1/(1-\lambda^{j-2})]^{n-1}\bfs_j$ for all $n \in {\mathbb N}$ and all $j \geq 3$.  If we choose $\bfu_1 = \bfe_k$ for some $k \geq 3$ then $\bfu_n = [1 - 1/(1-\lambda^{k-2})]^{n-1} \bfe_k$ and hence
$$
\| \bfu_n\|_{\infty}^{1/n} = [1/(1-\lambda^{k-2}) - 1]^{1-1/n} \rightarrow [1/(1-\lambda^{k-2}) - 1] = \lambda^{k-2}/(1 - \lambda^{k-2}) \leq \lambda/(1 - \lambda)
$$
as $n \rightarrow \infty$. It follows that $X_{\mbox{\scriptsize \rm reg}} \subset X$ is the closed subspace $c_0 \cap \mbox{sp}(\{ \bfe_3, \bfe_4, \ldots \})^a \subset c_0$.

Therefore $X = X_{\mbox{\scriptsize \rm sin}} \oplus X_{\mbox{\scriptsize \rm reg}} \cong X_{\mbox{\scriptsize \rm sin}} \times X_{\mbox{\scriptsize \rm reg}}$ and $Y = Y_{\mbox{\scriptsize \rm sin}} \oplus Y_{\mbox{\scriptsize \rm reg}} \cong Y_{\mbox{\scriptsize \rm sin}} \times Y_{\mbox{\scriptsize \rm reg}}$ where $Y_{\mbox{\scriptsize \rm sin}} = B_1(X_{\mbox{\scriptsize \rm sin}}) = \mbox{sp}(\{ \bfe_1, \bfe_2 \})$ and $Y_{\mbox{\scriptsize \rm reg}} = B_0(X_{\mbox{\scriptsize \rm reg}}) = c_0 \cap \mbox{sp}( \{ \bfe_3, \bfe_4, \ldots \})^a$. The corresponding key projections are $P = Q = [\bfe_1, \bfe_2, \bfzero, \bfzero, \ldots] \in {\mathcal B}(c_0)$ and the basic solution is found by solving $T_{-1}B_1 = P$ and $B_1T_{-1} = Q$ to find $T_{-1} \in {\mathcal B}(Y,X)$ and by solving $T_0B_0 = P^c$ and $B_0T_0 = Q^c$ to find $T_0 \in {\mathcal B}(Y,X)$. Elementary algebra gives
$$
T_{-1} = [-\bfe_1, \bfe_1 - \bfe_2, \bfzero, \bfzero, \ldots] \quad \mbox{and} \quad T_0 = [\bfzero,\bfzero, \{1/(1-\lambda)\}\bfe_3, \{1/(1-\lambda^2)\}\bfe_4,\ldots].
$$
The algebra is simply Gaussian elimination but the calculation of $T_0$ requires consideration of both $T_0B_0 = P^c$ and $B_0T_0 = Q^c$. Now we can calculate the entire Laurent series using $T_{-k} = (-1)^{k-1} (T_{-1}B_0)^{k-1}T_{-1} \in {\mathcal B}(Y,X)$ and $T_{\ell} = (-1)^{\ell} (T_0B_1)^{\ell} T_0 \in {\mathcal B}(Y,X)$ which gives
$$
T_{-2} = [\bfzero, \bfe_1, \bfzero, \bfzero,\ldots] \quad \mbox{and} \quad T_{-k} = \bigzero\ \mbox{for}\ k \geq 3
$$
and
$$
T_{\ell} = [ \bfzero, \bfzero, \{\lambda^{\ell}/(1-\lambda)^{\ell + 1}\}\bfe_3, \{\lambda^{2 \ell}/(1-\lambda^2)^{\ell+1}\}\bfe_4,\ldots]\ \mbox{for}\ \ell \in {\mathbb N}.
$$
Therefore $R(z)$ has a pole of order $2$ at $z=1$. The singular part is $R_{\mbox{\scriptsize \rm sin}}(z) = T_{-2}/(z-1)^2 + T_{-1}/(z-1)$ for $z \neq 1$ and the regular part is $R_{\mbox{\scriptsize \rm reg}}(z) = T_0 + T_1(z-1) + T_2(z-1)^2 + \cdots$ for $z \in D_{1/(1 - \lambda)}(1)$. Thus the Laurent series converges for $z \in D_{0,1/(1-\lambda)}(1)$.} $\hfill \Box$
\end{exa}

\section{An application to the control of electrical power generation}
\label{s:aepg}

We consider an application to the dynamic control of electrical power generation in a wind farm. The proposed feedback control is designed to mitigate the high volatility of power generation in wind farms observed by \cite{agr1}. The time series $\bfx$ is a non-autonomous $AR(1)$ process. Thus\textemdash strictly speaking\textemdash this example lies beyond the theoretical scope of the paper. However we will show that the new methods can be used to find $\bfx(t)$.

\begin{exa}[Power generation]
\label{ex:4}

{\rm Business Managers of Australian wind-farms sign short-term contracts with the Australian Energy Market Operator to supply an agreed number $H = H_{\mathcal T}$ of kilowatt-hours of electrical energy during the next trading interval $t \in {\mathcal T} = \{1,\ldots, N\}$. Let $\bfw = \{ \bfw(t)\}_{t \in {\mathbb N}-1} \in {\mathbb R}^{{\mathbb N}-1}$ be the number of kilowatt-hours of electrical energy collected on the interval $[0,t]$. To moderate the observed high volatility we seek a suitable feedback mechanism that encourages an appropriate target for the energy $\Delta \bfw(t) = \bfw(t) - \bfw(t-1)$ collected during each subinterval $[t-1,t]$ for $t \in {\mathcal T}$. Our aim will be to ensure that
$$
\bfx(t) = \bfw(t) - Ht/N \rightarrow 0
$$
as $t \rightarrow N$ where $\bfx(t)$ is the notional cumulative error at time $t \in {\mathcal T}$. Intuitively we could manage the error by setting $\bfx(t) - \bfx(t-1) = - p \bfx(t-1) / (t + q)$ for some $p, q \in {\mathbb N}$ with $p \leq q$ and all $t \in {\mathcal T}$. Thus we would set $\bfx(t) = [1 - p/ (t + q)] \bfx(t-1)$. Now $1 - p / (t + q) \in (0,1)$ for all $t \in {\mathcal T}$ and so either $\bfx(t-1) \leq \bfx(t) \leq 0$ or $0 \leq \bfx(t) \leq \bfx(t-1)$. Thus the notional error must decrease in magnitude. In practice we observe $\bfx_{\sbxi}(t) = \bfx(t) + \bxi(t)$ where $\bxi(t)$ is an inherent stochastic noise. Hence our target becomes
\begin{equation}
\label{e:(7.1)}
\bfx(t) = [ 1 - p/(t + q)] \bfx_{\sbxi}(t-1) \iff \bfx_{\sbxi}(t) = [1 - p/(t + q)] \bfx_{\sbxi}(t-1) + \bxi(t) 
\end{equation}
for all $t \in {\mathcal T}$. Suppose (\ref{e:(7.1)}) is true for all $t \in {\mathbb N}$ and assume $\bfx_{\sbxi}(0) = 0$. It follows that
$$
\mbox{$\sum_{t \in {\mathbb N}}$} \bfx_{\sbxi}(t)z^t = z\, \mbox{$\sum_{t \in {\mathbb N}}$} \bfx_{\sbxi}(t-1)z^{t-1} - (p/z^q) \left[ \mbox{$ \sum_{t \in {\mathbb N}}$} \bfx_{\sbxi}(t-1) z^{t+q}/(t+q) \right] + \mbox{$\sum_{t \in {\mathbb N}}$} \bxi(t)z^t
$$
which we can rewrite in terms of the ${\mathcal M}$-transform as
\begin{equation}
\label{e:(7.2)}
\bfX_{\sbxi}(z) = z \bfX_{\sbxi}(z) - W \bfX_{\sbxi}(z) + \bXi(z) - \bxi(0) \iff [(1 - z) I + W] \bfX_{\sbxi}(z) = \bXi(z) - \bxi(0)
\end{equation}
where $W$ is a weighted Volterra operator defined by the formula $W \bfF(z) = (p/z^q) \int_0^z F(\zeta)\zeta^q d\zeta$ for all $\bfF(z) = \sum_{t \in {\mathbb N}-1} \bff(t)z^t$ which are analytic when $z \in D_r(0)$ for some $r > 0$. An intuitive application of the Neumann expansion suggests that we may have
\begin{equation}
\label{e:(7.3)}
[I(z-1) + W]^{-1} = \left\{ I - (1-z)^{-1}W + [(1-z)^{-1}W]^2 - [(1-z)^{-1}W]^3 + \cdots \right\} (1 - z)^{-1}
\end{equation}
for $z \neq 1$. This is not entirely true but a detailed justification of (\ref{e:(7.3)}) for $z \in D_1(0)$ is given in Section~\ref{s:wvo} of the Appendix. Order is important because $(1 - z)^{-1} W \neq W (1-z)^{-1}$. The formula (\ref{e:(7.3)}) suggests that $[I(z-1) + W]^{-1}$ has an essential singularity at $z = 1$ and that $\bfx_{\sbxi} \sim I(\infty)$. If $z \in D_1(0)$ and we write $(1-z)^{-1} = 1+z+z^2+\cdots$ then straightforward calculations give
$$
(1 - z)^{-1}[ \bXi(z) - \bxi(0) ] = \bxi(1)z + [\bxi(1) + \bxi(2)]z^2 + [\bxi(1) + \bxi(2) + \bxi(3)]z^3 + \cdots,
$$
\begin{eqnarray*}
\lefteqn{(1-z)^{-1}W(1-z)^{-1}[\bXi(z) - \bxi(0)] } \\
& = & \{p\bxi(1)/(q+2)\}z^2 + \{p\bxi(1)/(q+2) + p[\bxi(1) + \bxi(2)]/(q+3)\}z^3 + \cdots,
\end{eqnarray*}
\begin{eqnarray*}
\lefteqn{[(1-z)^{-1}W]^2(1-z)^{-1}[\bXi(z) - \bxi(0)] = \{p^2\bxi(1)/[(q+2)(q+3)]\}z^3 } \\
& & + \{p^2\bxi(1)/[(q+2)(q+3)] + \left[\{p\bxi(1)/(q+2) + p[\bxi(1) + \bxi(2)]/(q+3)\}/(q+4) \right] \} z^4 + \cdots
\end{eqnarray*}
and so on. By collecting the coefficients of the various powers we have
\begin{eqnarray*}
\lefteqn{[(1-z)I + W]^{-1}[\bXi(z) - \bxi(0)]} \\
& = & \bxi(1)z + \{[1 - p/(q+2)]\bxi(1) + \bxi(2)\}z^2 \\ 
& & +\{ [1 - p/(q+3)][1 - p/(q+2)]\bxi(1) + [1 - p/(q+3)]\bxi(2) + \bxi(3)\} z^3 + \cdots \\
& = & \bfx_{\sbxi}(1)z + \bfx_{\sbxi}(2)z^2 + \bfx_{\sbxi}(3)z^3 + \cdots.
\end{eqnarray*}
In general we can show that
$$
\bfx_{\sbxi}(t) = \left\{ \begin{array}{ll}
\bxi(1) & \mbox{for}\ t=1 \\
\bxi(t) + \mbox{$\sum_{s=1}^{t-1}$} \left\{\rule{0cm}{0.4cm} \mbox{$\prod_{r=t-s+1}^{t}$} [1 - p/(q+r)] \right\} \bxi(t-s) & \mbox{for}\ t \in {\mathbb N}+1. \end{array} \right.
$$
Therefore $\bfx_{\sbxi}(t)$ may grow without bound as $t \rightarrow \infty$. It can be shown that $\Delta^k \bfx_{\sbxi}(t)$ may also grow without bound for all $k \in {\mathbb N}$ as $t \rightarrow \infty$. Thus $\bfx_{\sbxi} \sim I(\infty)$.  The corresponding increment is $\Delta \bfw_{\sbxi}(t) = - p\bfx_{\sbxi}(t-1)/(t+q) + H/N$ for all $t \in {\mathcal T}$. In practice $p, q$ could be adjusted to get the best results. The detailed calculations are justified in Section~\ref{s:wvo} of the Appendix.} $\hfill \Box$
\end{exa} 

\begin{rem}
\label{r:5}
In retrospect we could change the dependent variable to $\bfy_{\sbxi}(t) = (t+q)^{(q-p)}\bfx_{\sbxi}(t)$ in Example~\ref{ex:4}.  Now $\bfy_{\sbxi}(t) - \bfy_{\sbxi}(t-1) = (t+q)^{(q-p)}\bxi(t)$ for all $t \in {\mathbb N}$ with $\bfy_{\sbxi}(0) = 0$ where the factorial power is defined by $(t+r)^{(s)} = (t+r)(t+r-1)\cdots (t+r-s+1)$ for all $t ,r, s \in {\mathbb N-1}$ with $s \in [r, t+r]$. Hence $\bfy_{\sbxi}$ is an autonomous $AR(1)$ process with an unbounded $I(\infty)$ noise term. A simple recursion gives $\bfy_{\sbxi}(t) = (t+q)^{(q-p)}\bxi(t) + (t-1+q)^{(q-p)}\bxi(t-1) + \cdots (1+q)^{(q-p)}\bxi(1)$ for all $t \in {\mathbb N}$.  Hence $\bfy_{\sbxi}(t)$ may grow without bound as $t \rightarrow \infty$.  It can be shown that $\Delta^k \bfy_{\sbxi}(t)$ may also grow without bound for each $k \in {\mathbb N}$ as $t \rightarrow \infty$. Therefore $\bfy_{\sbxi} \sim I(\infty)$.
\end{rem}

\section{Conclusions}
\label{s:c}

We have used the new formula $R(z) = [(z-1)I_X + T_{-1}B_0]^{-1}T_{-1} + [I_X + T_0B_1(z-1)]^{-1}T_0$ for the resolvent of the autoregressive polynomial to establish an extended GJRT that applies to all integrated autoregressive processes irrespective of the dimension of the space, the number of stochastic trends and cointegrating relations, and the order of integration, all of which may be finite or infinite. The present results are a substantial extension of those reported in \cite{bea1}, \cite{bea2}, \cite{fra5}, and \cite{seo1}.  Standard linear algebra can be used to find the key spectral separation projections in finite-dimensional problems but it may be necessary to use infinite-length Jordan chains in more general problems.

\section{Acknowledgements}
\label{s:ack}

Massimo Franchi gratefully acknowledges partial financial support from Ministero dell' Universit\`{a} e della Ricerca for the project PRIN 20223725WE \textit{Methodological and computational issues in large-scale time series models for economics and finance}. Phil Howlett and Brendan Beare would like to thank the organisers of ANZESG 2023 and the participants for critiquing the initial presentation of this research.

\appendix

\section{Appendix: Notation and terminology}
\label{s:nt}

We use the following notation throughout. In general $X, Y$ will denote complex Banach spaces and $(\Omega, \Sigma, \mu)$ will denote a probability space.

\begin{itemize}

\item ${\mathbb N} = \{ 1,2,\ldots\}$ is the set of natural numbers, ${\mathbb N} - 1$ is the set of nonnegative integers, ${\mathbb Z} =(-{\mathbb N}) \cup ({\mathbb N}-1)$ is the set of integers, ${\mathbb N} + m$ is the set $\{1+m, 2+m, \ldots \}$ for each $m \in {\mathbb Z}$, ${\mathbb R}$ is the set of real numbers, and ${\mathbb C}$ is the set of complex numbers. If $z \in {\mathbb C}$ then $\overline{z} \in {\mathbb C}$ denotes the complex conjugate.

\item If $D \subset {\mathbb C}$ we follow \citet[p.~3]{yos1} and write $D^a$ to denote the closure of $D$.

\item For $z_0 \in {\mathbb C}$ and $s, r \in [0, \infty] \subset {\mathbb R} \cup \{ \infty\}$ with $s < r$ we write $D_{\,r}(z_0) = \{ z \in {\mathbb C} \mid |z-z_0| < r\}$ and $D_{\,s,r}(z_0) = \{ z \in {\mathbb C} \mid s < |z-z_0| <r\}$.

\item A linear operator $A: X \rightarrow Y$ is bounded if $\|A \bfx\|_Y \leq c \| \bfx \|_X$ for some $c \in {\mathbb R}$ with $c > 0$. The set of all such bounded linear operators is denoted by ${\mathcal B}(X,Y)$.  If $A \in {\mathcal B}(X,Y)$ and $T \subset Y$ we write $A^{-1}(T) = \{ \bfx \in X \mid A\bfx \in T\}$ whether or not $A$ is invertible.  The set ${\mathcal B}(X,X)$ is denoted by ${\mathcal B}(X)$.  The identity operator $I_X \in {\mathcal B}(X)$ is defined by $I_X \bfx = \bfx$ for all $\bfx \in X$. If $P \in {\mathcal B}(X)$ is a projection we write $P^c = I_X - P$ to denote the complementary projection. We may write $I$ instead of $I_X$ if $X$ is clear from the context.  

\item We follow \citet[Sections 5.6, 5.7, pp.~115\textendash 118]{lue1} and use the scalar product notation $\langle \bfx, \bff \rangle \in {\mathbb C}$ to denote the value of the bounded linear functional $\bff \in X^*$ at the point $\bfx \in X$. If $S \subset X$ is a subset of $X$ then $S^{\perp} \subset X^*$ is the subspace defined by $S^{\perp} = \{ \bff \in X^*\mid \langle \bfs, \bff \rangle = 0\ \mbox{for all}\ \bfs \in S\}$.

\item A stochastic time series on $X$ is a sequence of strongly measurable maps from $\Omega$ to $X$. The stochastic time series $\bfx = \{\bfx(t)\}_{t \in {\mathbb Z}} :\Omega \rightarrow X^{\mathbb Z}$ is a series with value $\bfx(\bomega) = \{ \bfx(t) \}_{t \in {\mathbb Z}}(\bomega) = \{\bfx(\bomega,t)\}_{t \in {\mathbb Z}} \in X^{\mathbb Z}$ that depends on the outcome $\bomega \in \Omega$ of some random process in the probability space $(\Omega, \Sigma, \mu)$. Thus the value of the series is simply defined as the corresponding series of values. If $A \in {\mathcal B}(X,Y)$ is a bounded linear operator then the stochastic time series $\bfy = A\bfx$ is defined by $\bfy(\bomega) = [A \{\bfx(t)\}_{t \in {\mathbb Z}}](\bomega) = \{A\bfx(\bomega,t)\}_{t \in {\mathbb Z}} \in Y^{\mathbb Z}$ for each $\bomega \in \Omega$. The dependence on $\bomega$ is normally suppressed unless it is directly relevant to the issue at hand.

\item According to \citet[p~94]{doo1} a stochastic time series $\bfx = \{ \bfx(t) \}_{t \in {\mathbb Z}} \in X^{\mathbb Z}$ is {\em strictly stationary} if the multivariate joint distribution of the variables $\bfx(t_1+q),\ldots,\bfx(t_p+q)$ for each finite collection of times $\{t_1, \ldots, t_p\} \in {\mathbb Z}^p$ does not depend on $q \in {\mathbb Z}$.

\item For each stochastic time series $\bfu \in X^{\mathbb Z}$ we define the associated observable time series $\bfu_+ = \{ \bfu_+(t) \}_{t \in {\mathbb Z}} \in X^{\mathbb Z}$ by setting $\bfu_+(t) = \bfu(t)$ for $t \in {\mathbb N}-1$ and $\bfu_+(t) = \bfzero$ otherwise.

\item For each $\bfu \in X^{\mathbb Z}$ the time series $\Delta \bfu \in X^{\mathbb Z}$ of first differences is defined using the (backward) difference operator $\Delta \bfu(t) = \bfu(t)-\bfu(t-1)$ and the corresponding series delayed by one unit of time $L\bfu \in X^{\mathbb Z}$ is defined using the lag operator $L \bfu(t) = \bfu(t-1)$ for each $t \in {\mathbb Z}$. The operator identity $\Delta = I_X - L$ shows that
\begin{equation}
\label{e:(A.1)}
\Delta^{-k} \bfu(t) = (I_X - L)^{-k} \bfu(t) = \mbox{$\sum_{s \in {\mathbb N}-1}$} \mbox{$\binom{k+s-1}{s}$} L^s \bfu(t) = \mbox{$\sum_{s \in {\mathbb N}-1}$} \mbox{$\binom{k+s-1}{s}$} \bfu(t-s)
\end{equation}
for all $k \in {\mathbb N}$, provided ${\mathbb E}[ \| \bfu(t-s) \| ] \leq c(t) /r^s$ for some $c(t), r > 0$ and all sufficiently large $s \in {\mathbb N}$ and that
\begin{equation}
\label{e:(A.2)}
\Delta^{\ell} \bfu(t) = (I_X - L)^{\ell}\bfu(t) = \mbox{$\sum_{s=0}^{\ell}$} \mbox{$\binom{\ell}{s}$} (-1)^s L^s \bfu(t) = \mbox{$\sum_{s=0}^{\ell}$} \mbox{$\binom{\ell}{s}$} (-1)^s \bfu(t-s)
\end{equation}
for all $\ell \in {\mathbb N}-1$. In each case $\binom{p}{q}$ is the binomial coefficient for $p,q \in {\mathbb N}-1$ with $p \geq q$.

\item A stochastic time series $\bxi = \{\bxi(\bomega,t)\}_{t \in {\mathbb Z}} \in Y^{\mathbb Z}$ is called a {\em strong white noise} process if it is an independent identically distributed (i.i.d.)~process with $0<\sigma^2 = {\mathbb E} \left[ \| \bxi \|^2 \right] = \mbox{$\int_{\Omega}$} \| \bxi(\bomega,t)\|^2 d \mu(\bomega) < \infty$ and $\bmu = {\mathbb E} \left[ \bxi \right] = \mbox{$\int_{\Omega}$} \bxi(\bomega,t) d \mu(\bomega) = \bfzero$ for each $t \in {\mathbb Z}$, where the latter integral is a Bochner integral \citep[Chapter V, Sections 4 and 5, pp. 130\textendash 134]{yos1}. The noise process is not necessarily Gaussian.  A strong white noise process is strictly stationary.

\item A stochastic time series $\bfx=\{\bfx(\bomega, t)\}_{t \in {\mathbb Z}} \in X^{\mathbb Z}$ is said to be {\em integrated of order zero}, written $\bfx \sim I(0)$, if there exists a strong white noise process $\{\bxi(\bomega, t)\}_{t \in \mathbb Z} \in Y^{\mathbb Z}$ and a sequence of bounded linear operators $\{A_s\}_{s\in\mathbb N-1}\in {\mathcal B}(Y,X)^{\mathbb N-1}$ with $\sum_{s \in {\mathbb N}-1} \left \Vert A_s \right \Vert < \infty$ and $\sum_{s\in\mathbb N-1}A_s \neq \bigzero$ so that $\bfx(\bomega,t)=\sum_{s\in\mathbb N-1} A_s \bxi(\bomega, t-s)$ for all $t \in\mathbb Z$.  In \citet[p 183]{bos1} a series in this form is called a {\em standard linear process}.  A standard linear process is strictly stationary and converges almost surely.  See \citet[Lemma 7.1 pp 182\textendash 183]{bos1}.

\item A stochastic linear process $\bfx$ is said to be {\em integrated of order $d$} for some $d \in {\mathbb N}$, written $\bfx \sim I(d)$, if $d$ is the smallest natural number such that $\Delta^d \bfx \sim I(0)$.  A stochastic linear process $\bfx \in X^{\mathbb Z}$ is said to be {\em integrated of order infinity}, written $\bfx \sim I(\infty)$ if there is no $d \in {\mathbb N}$ such that $\Delta^d \bfx \sim I(0)$. 

\end{itemize}

\section{Appendix: Proofs of the main results}
\label{s:pmr}

{\bf Proof of Lemma~\ref{lem:rcon}.} Suppose $R(z)$ is analytic on $D_{1 + \epsilon}(0) \setminus \{1\}$ for some $\epsilon > 0$. Therefore $R(z)$ is analytic on $[\, D_{\, 1}(0)^a \setminus \{1\}\,] \cup D_{\, 0,\, \epsilon}(1) \subset D_{1 + \epsilon}(0) \setminus \{1\}$. Now suppose $R(z)$ is analytic on $[\, D_{\, 1}(0)^a \setminus \{1\}\,] \cup D_{\, 0,\, \delta}(1)$ for some $\delta > 0$. Define the compact set ${\mathcal C} = \{ z \in {\mathbb C} \mid |z| = 1\ \mbox{and}\ z \notin D_{\, \delta/2}(1) \}$. The resolvent $R(z)$ is analytic on ${\mathcal C}$ and so for each $w \in {\mathcal C}$ we can find $\epsilon(w) > 0$ such that $R(z)$ is analytic on the open neighbourhood $D_{\, \epsilon(w)}(w)$. Therefore ${\mathcal C} \subset \mbox{$\bigcup_{w \in {\mathcal C}}$}\, D_{\, \epsilon(w)}(w)$.  The set ${\mathcal C}$ is compact so there exists a finite subcollection of open neighbourhoods $\{ D_{\epsilon_i}(w_i) \}_{i=1}^n$ with $w_i \in {\mathcal C}$ and $\epsilon_i > 0$ for each $i=1,\ldots,n$ such that ${\mathcal C} \subset \mbox{$\bigcup_{i=1}^n$}\, D_{\, \epsilon_i}(w_i)$.  The set $\mbox{$\bigcup_{i=1}^n$}\, D_{\, \epsilon_i}(w_i)$ is open so there is some $\epsilon \in (0, \delta/2)$ such that $N({\mathcal C}, \epsilon) = \{ z \in {\mathbb C} \mid |z - w| < \epsilon\ \mbox{for some}\ w \in {\mathcal C}\} \subset \mbox{$\bigcup_{i=1}^n$}\, D_{\, \epsilon_i}(w_i)$. Thus $R(z)$ is analytic for $z \in N({\mathcal C}, \epsilon)$. Now $D_{1 + \epsilon}(0) \setminus \{1\} \subset N({\mathcal C}, \epsilon) \cup [\, D_{\, 1}(0)^a \setminus \{1\}\,] \cup D_{\, 0,\, \delta}(1)$. Consequently $R(z)$ is analytic for $D_{1 + \epsilon}(0) \setminus \{1\}$. $\hfill \Box$.

\vspace{0.5cm}

{\bf Proof of Lemma~\ref{lem:mtrr}.} We observe that $PW(z) = \sum_{t \in {\mathbb N}-1} PW_t z^t$ for all $z \in D_{1+\epsilon}(0)$. However we also know that $PT_{\ell} = \bigzero$ for all $\ell \in {\mathbb N}-1$ and so $PW(z) = PR_{\mbox{\scriptsize \rm reg}}(z) = \mbox{$\sum_{\ell \in {\mathbb N}-1}$} PT_{\ell}(z-1)^{\ell} = \bigzero$ for all $z \in D_{\, \epsilon}(1)$. Since $PW(z) \equiv \bigzero$ on a nontrivial open set it follows by analytic continuation that $PW(z) \equiv \bigzero$ for all $z \in D_{1+\epsilon}(0)$ and hence $PW_t = \bigzero$ for all $t \in {\mathbb N}-1$.  A similar argument shows that $W(z)Q = \bigzero$ for all $z \in D_{1+\epsilon}(0)$. Therefore $W_tQ = \bigzero$ for all $t \in {\mathbb N}-1$. $\hfill \Box$

\vspace{0.5cm}

{\bf Proof of Proposition~\ref{p:1}.} Consider the term $\bfX_{\mbox{\scriptsize \rm sin}}(z)$. It follows from (\ref{e:(3.11)}) that
\begin{eqnarray}
\label{e:(B.1)}
\lefteqn{ \bfX_{\mbox{\scriptsize \rm sin}}(z) = R_{\mbox{\scriptsize \rm sin}}(z) [\bXi(z) - B_1 \bfc] = \mbox{$\sum_{k \in {\mathbb N}}$} T_{-k}(z - 1)^{-k} \left[ \bXi(z) - B_1 \bfc \right] } \hspace{4.53cm} \nonumber \\
& & = \mbox{$\sum_{k \in {\mathbb N}}$} T_{-k}\, (-1)^k {\mathcal M} \left[ \{ \Delta^{-k} (\bxi - B_1 \bfc\, \bdelta)_+(t) \}_{t \in {\mathbb N}-1} \right](z). \hspace{1cm}
\end{eqnarray} 
Now we can apply (\ref{e:(3.4)}), (\ref{e:(3.9)}) and (\ref{e:(4.4)}) to deduce that
\begin{eqnarray}
\label{e:(B.2)}
\bfX_{\mbox{\scriptsize \rm sin}}(z) & = & {\mathcal M} \left[ \{ \mbox{$\sum_{k \in {\mathbb N}}$} (-1)^kT_{-k} \Delta^{-k} \bxi_+(t) \}_{t \in {\mathbb N}-1} \right](z) \nonumber \\
& & \hspace{2cm} - {\mathcal M} \left[ \{ \mbox{$\sum_{k \in {\mathbb N}} \binom{t+ k-1}{k-1}$}(T_{-1}B_0)^{k-1}T_{-1} B_1 \bfc \}_{t \in {\mathbb N}-1} \right](z) \nonumber \\
& = & {\mathcal M} \left[ \mbox{$\sum_{k \in {\mathbb N}}$} (-1)^k T_{-k} \Delta^{-k} \bxi_+(t) - U_tB_1 \bfc \right] (z). 
\end{eqnarray}
We justify convergence of (\ref{e:(B.2)}) as follows. For each $\eta > 0$, we can find $c_{\eta} > 0$ such that $
\|T_{-k} \| \leq c_{\eta} \eta^k$ for all $ k \in {\mathbb N}$. If we choose $\eta \in (0,1)$ then
\begin{eqnarray*}
{\mathbb E} \left[ \| \mbox{$\sum_{k \in {\mathbb N}}$} (-1)^k T_{-k} \Delta^{-k} \bxi_+(t) \| \right] & = & {\mathbb E} \left[ \, \| \mbox{$\sum_{k \in {\mathbb N}}$} (-1)^k T_{-k} \mbox{$\sum_{s=0}^t \binom{s + k-1}{k-1}$} \bxi(t-s) \| \right] \\
& \leq & {\mathbb E} \left[ \| \bxi \| \right] \cdot \mbox{$\sum_{k \in {\mathbb N}}$}\, c_{\eta} \eta^k \cdot \mbox{$\sum_{s=0}^t \binom{s + k-1}{k-1}$}  \\
& \leq & c_{\eta}\, {\mathbb E} \left[ \| \bxi \| \right] \cdot \mbox{$\sum_{k \in {\mathbb N}} \binom{t+ k}{k}$} \eta^k \\
& = & c_{\eta}\, {\mathbb E} \left[ \| \bxi \| \right] \cdot \left[ 1 / (1 - \eta)^{t+1} - 1 \right]
\end{eqnarray*}
is well defined and finite. Therefore $\bfX_{\mbox{\scriptsize \rm sin}}(z) = {\mathcal M} \left[ \{ \bfx_{\mbox{\scriptsize \rm sin}}(t) \}_{t \in {\mathbb N}-1} \right](z)$ where
\begin{equation}
\label{e:(B.3)}
\bfx_{\mbox{\scriptsize \rm sin}}(t) = \mbox{$\sum_{k \in {\mathbb N}}$} (-1)^k T_{-k} \Delta^{-k} \bxi_+(t) - U_t B_1 \bfc \in P(X)
\end{equation}
for each $t \in {\mathbb N}-1$. Now consider the term $\bfX_{\mbox{\scriptsize \rm reg}}(z)$. We have
\begin{eqnarray}
\label{e:(B.4)}
\bfX_{\mbox{\scriptsize \rm reg}}(z) & = & W(z)[ G(z) - B_1 \bfc] \nonumber \\
& = & \mbox{$\sum_{s \in {\mathbb N}-1}$} W_s z^s\, \mbox{$\sum_{r \in {\mathbb N}-1}$} (\bxi - B_1\bfc\, \bdelta)_+(r) z^r \nonumber \\
& = & \mbox{$\sum_{t \in {\mathbb N}-1}$} \left[ \mbox{$\sum_{s \in {\mathbb N}-1}$} W_s\, (\bxi - B_1 \bfc\, \bdelta)_+(t-s) \right] z^t \nonumber \\
& = & \mbox{$\sum_{t \in {\mathbb N}-1}$} \left[ \mbox{$\sum_{s=0}^t$} W_s\, \bxi(t-s) - W_t\, B_1 \bfc \right] z^t
\end{eqnarray}
for all $z \in D_{1 + \epsilon}(0)$. Therefore $\bfX_{\mbox{\scriptsize \rm reg}}(z) = {\mathcal M} \left[ \bfx_{\mbox{\scriptsize \rm reg}}(t) \right] (z)$ where
\begin{equation}
\label{e:(B.5)}
\bfx_{\mbox{\scriptsize \rm reg}}(t) = \mbox{$\sum_{s=0}^t$}\, W_s\, \bxi(t-s) - W_t\, B_1 \bfc \in P^c(X)
\end{equation}
for all $t \in {\mathbb N}-1$.  We know that $\|W_s\| < c/(1+\epsilon)^s$ and so convergence is straightforward with
$$
{\mathbb E}[ \| \mbox{$\sum_{s=0}^t$}\, W_s\, \bxi(t-s) \|] \leq {\mathbb E}[ \| \bxi\|] \mbox{$\sum_{s=0}^t$} c/(1 + \epsilon)^s \leq c\, {\mathbb E}[ \| \bxi\|] (1 + \epsilon)/\epsilon
$$ 
for all $t \in {\mathbb N}-1$.  Finally we set $\bfx(t) = \bfx_{\mbox{\scriptsize \rm sin}}(t) + \bfx_{\mbox{\scriptsize \rm reg}}(t)$ for all $t \in {\mathbb N}-1$.   $\hfill \Box$

\vspace{0.5cm}

{\bf Proof of Corollary~\ref{c:1}.} We have $R_{\mbox{\scriptsize \rm sin}}(z) = [(z-1)I_X + T_{-1}B_0]^{-1}T_{-1} = \mbox{$\sum_{t \in {\mathbb N}-1}$} U_tz^t$ for $z \in D_1(0)$ where $PU_t = U_t$ for all $t \in {\mathbb N}-1$. Therefore $PR_{\mbox{\scriptsize \rm sin}}(z) = R_{\mbox{\scriptsize \rm sin}}(z)$.  The equation
\begin{equation}
\label{e:(B.6)}
\bfX_{\mbox{\scriptsize \rm sin}}(z) = R_{\mbox{\scriptsize \rm sin}}(z)[\bXi(z) - B_1 \bfc]
\end{equation}
shows that we also have $P\bfX_{\mbox{\scriptsize \rm sin}}(z) = \bfX_{\mbox{\scriptsize \rm sin}}(z)$.  If we multiply both sides of (\ref{e:(B.6)}) on the left by $B_1[(z-1)I_X+T_{-1}B_0]$ and use the identity $T_0B_0 + T_{-1}B_1 = I_X$ we get
$$
B_1[ (z-1)I_X + T_{-1}B_1][T_0B_0 + T_{-1}B_1]\bfX_{\mbox{\scriptsize \rm sin}}(z) = B_1T_{-1} \bXi(z) - B_1T_{-1}\bfc
$$
for all $z \in D_1(0)$.  Now we can use the relationships $B_0 = A_0 + A_1$ and $B_1 = A_1$, the identity $T_{-1}B_0T_0 = \bigzero$ and the definitions $P = T_{-1}B_1$ and $Q = B_1T_{-1}$ to deduce that
$$
[QA_0P + QA_1P z] \bfX_{\mbox{\scriptsize \rm sin}}(z) = Q\bXi(z) - QB_1\bfc
$$
for all $z \in D_1(0)$.  Finally we have $QB_1 = B_1T_{-1}[B_1T_{-1} + B_0T_0]B_1 = B_1 P^2  = A_1 P$ from which it follows that $QB_1 = Q^2B_1 = QA_1P$.  Thus we obtain
\begin{equation}
\label{e:(B.7)}
[QA_0P + QA_1P z] \bfX_{\mbox{\scriptsize \rm sin}}(z) = Q\bXi(z) - QA_1P \bfc
\end{equation}
for all $z \in D_1(0)$.  This means that $\bfx_{\mbox{\scriptsize \rm sin}}(t) \in P(X)$ satisfies the $AR(1)$ equation
\begin{equation}
\label{e:(B.8)}
QA_0P \bfx_{\mbox{\scriptsize \rm sin}}(t)  +  QA_1P \bfx_{\mbox{\scriptsize \rm sin}}(t - 1) = Q \bxi(t)
\end{equation}
for all $t \in {\mathbb N} - 1$ with $\bfx_{\mbox{\scriptsize \rm sin}}(-1) = QA_1P\bfc$.  We know that $QB_iP + Q^cB_iP^c = B_i$ for each $i=0,1$.  It follows that $QA_iP + Q^cA_iP^c = A_i$ for each $i=0,1$.  Therefore
$$
A_0\bfx(t) + A_1 \bfx(t-1) - [QA_0P \bfx_{\mbox{\scriptsize \rm sin}}(t)(t) + QA_1P \bfx_{\mbox{\scriptsize \rm sin}}(t-1) = \bxi(t) - Q\bxi(t)
$$
for all $t \in {\mathbb N}- 1$.  We have $A_i \bfx(t) = (QA_iP + Q^cA_iP^c)[\bfx_{\mbox{\scriptsize \rm sin}}(t) + \bfx_{\mbox{\scriptsize \rm reg}}(t)] = QA_iP \bfx_{\mbox{\scriptsize \rm sin}}(t) + Q^cA_iP^c \bfx_{\mbox{\scriptsize \rm reg}}(t)$ for each $i=0,1$.  Therefore the previous equation becomes
\begin{equation}
\label{e:(B.9)}
Q^cA_0P^c \bfx_{\mbox{\scriptsize \rm reg}}(t)  +  Q^cA_1P^c \bfx_{\mbox{\scriptsize \rm reg}}(t - 1) = Q^c \bxi(t)
\end{equation}
for all $t \in {\mathbb N}-1$ with $\bfx_{\mbox{\scriptsize \rm reg}}(-1) = \bfc - Q\bfc = Q^c \bfc$.  If $\bfx \in X^{{\mathbb N}-1}$ and Assumption~\ref{a:1} is true we have shown that $\bfx = \bfx_{\mbox{\scriptsize \rm sin}} + \bfx_{\mbox{\scriptsize \rm reg}}$ where $\bfx_{\mbox{\scriptsize \rm sin}} \in P(X)^{{\mathbb N}-1}$ and $\bfx_{\mbox{\scriptsize \rm reg}} \in P^c(X)^{{\mathbb N}-1}$.  Similar algebraic manipulations show that the reverse implication is also true.

\vspace{0.5cm}

{\bf Proof of Corollary~\ref{c:2}.}  By recursive application of (\ref{e:(4.16)}) we obtain
\begin{equation}
\label{e:(B.10)}
\bfx_{\mbox{\scriptsize \rm reg}}(t) = \Phi^{t+r+1} \bfx_{\mbox{\scriptsize \rm reg}}(-1-r) + \mbox{$\sum_{s=0}^{t+r}$} \Phi^s \bfeta(t-s)
\end{equation}
for all $r \in {\mathbb N}-1$.  In particular for $t = -1$ we have
$$
P^c \bfc = \Phi^r \bfx_{\mbox{\scriptsize \rm reg}}(-1-r) + \mbox{$\sum_{s=0}^{r-1}$} \Phi^s \bfeta(-1-s)
$$
for all $r \in {\mathbb N}-1$.  We know that
$$
\mathbb E \left[ \Vert \mbox{$\sum_{s=0}^{r-1}$} \Phi^s \bfeta(-1-s) \Vert \right] \leq \mbox{$\sum_{s=0}^{r-1}$}\,k /(1 + \epsilon)^s \leq k(1 + 1/\epsilon)
$$
for some $k > 0$ and all $r \in {\mathbb N}$.  It follows that $\sum_{s=0}^{r-1} \Phi^s \bfeta(-1-s)$ converges almost surely as $r \rightarrow \infty$.  Therefore $\bfz_r = \Phi^r \bfx_{\mbox{\scriptsize \rm reg}}(-1-r)$ also converges almost surely as $r \rightarrow \infty$.  Thus we obtain
$$
\bfz_{\infty} = \lim_{r \rightarrow \infty} \bfz_r = \lim_{r \rightarrow \infty} \Phi^r \bfx_{\mbox{\scriptsize \rm reg}}(-1-r) = P^c \bfc - \mbox{$\sum_{s \in {\mathbb N}-1}$} \Phi^s \bfeta(-1-s) \in P^c(X)
$$
If we return to (\ref{e:(B.10)}) we now have
\begin{equation}
\label{e:(B.11)}
\bfx_{\mbox{\scriptsize \rm reg}}(t) = \Phi^{t+1} \cdot \Phi^r \bfx_{\mbox{\scriptsize \rm reg}}(-1-r) + \mbox{$\sum_{s=0}^{t+r}$} \Phi^s \bfeta(t-s)
\end{equation}
for all $r \in {\mathbb N}-1$.  By taking the limit in (\ref{e:(B.11)}) as $r \rightarrow \infty$ we get
$$
\bfx_{\mbox{\scriptsize \rm reg}}(t) = \Phi^{t+1} \bfz_{\infty} + \mbox{$\sum_{s \in {\mathbb N}-1}$} \Phi^s \bfeta(t-s)
$$
for all $t \in {\mathbb Z}$. If $R(z)$ has a pole of order $d$ at $z = 1$ then
\begin{eqnarray*}
U_t & = & (-1)(I_X - T_{-1}B_0)^{-t-1} \\
& = & \mbox{$\sum_{k \in {\mathbb N}-1} \binom{t+k}{k}$} (T_{-1}B_0)^k T_{-1} \\
& = & \mbox{$\sum_{k = 0}^{d-1} \binom{t+k}{k}$} (-1)^k T_{-k-1} 
\end{eqnarray*}
because $T_{-k-1} = 0$ for $k > d-1$.  If $\bfz_{\infty} = \bfzero$ it follows from (\ref{e:(4.13)}) that
$$
\bfx(t) =  \mbox{$\sum_{k = 0}^d$} (-1)^k T_{-k} \Delta^{-k} \bxi_+(t) + \mbox{$\sum_{k = 0}^{d-1}$} (-1)^{k+1} T_{-k-1} \mbox{$ \binom{t+k}{k}$} B_1 \bfc + \mbox{$\sum_{s \in {\mathbb N}-1}$} \Phi^s \bfeta(t-s) 
$$
for all $t \in {\mathbb N}-1$.  Now $\Delta \binom{t+k}{k} = \binom{t+k}{k} - \binom{t-1+k}{k} = \binom{t + k - 1}{k-1}$.  An inductive argument shows that $\Delta^k \binom{t+k}{k} = 1$.  Therefore  $\Delta^d \binom{t+k}{k} = 0$ for all $k < d$.  Some elementary algebra gives 
$$
\Delta^d \bfx(t) = T_{-d} \bxi(t) +  \mbox{$\sum_{s \in {\mathbb N}-1}$} \Phi^s \Delta^d \bfeta(t-s)
$$
for all $t \in {\mathbb N}-1+d$.  It follows that $\bfx \sim I(d)$, $\bfx_{\mbox{\scriptsize \rm sin}} \sim I(d)$ and $\bfx_{\mbox{\scriptsize \rm reg}} \sim I(0)$.  If $\bfz_{\infty} = \bfzero$ but $R(z)$ has an isolated essential singularity at $z = 1$ then no finite number of differences will reduce $\bfx_{\mbox{\scriptsize \rm sin}}(t)$ to a strictly stationary standard linear process.  Therefore $\bfx \sim I(\infty)$, $\bfx_{\mbox{\scriptsize \rm sin}} \sim I(\infty)$ and $\bfx_{\mbox{\scriptsize \rm reg}} \sim I(0)$.  $\hfill \Box$

\section{Appendix: The weighted Volterra operator}
\label{s:wvo}

The weighted Volterra operator $W$ is defined by $W[\bfF(z)] = (p/z^q) \mbox{$\int_0^z$} \bfF(\zeta) \zeta^q d\zeta$ for all functions $\bfF(z) = \bff(0) + \bff(1)z + \bff(2)z^2 + \cdots$ which are analytic on some region $D_r(0)$ where $r \in (0,1]$.  We wish to justify the formula $S = [(1-z)I + W]^{-1}$ where $S$ is defined by
\begin{equation}
\label{e:(C.1)}
S\bfG(z) \coloneqq \left[(1-z)^{-1} - [(1-z)^{-1}W](1-z)^{-1} + [(1-z)^{-1}W]^2(1-z)^{-1} - \cdots \right]\bfG(z)
\end{equation}
for all $\bfG(z) = \bfg(0) + \bfg(1)z + \bfg(2)z^2 + \cdots$ which are analytic for $z \in D_r(0)$.  Consider the operator $[(1-z)I + W]$.  Direct calculations show that $[(1-z)I + W] \bfF(z) = \bfG(z)$ where 
\begin{equation}
\label{e:(C.2)}
\bfg(0) = \bff(0), \ \bfg(1) = \bff(1) - [1 - p/(q+1)]\bff(0), \ \bfg(2) = \bff(2) - [1 - p/(q+2)]\bff(1), \cdots
\end{equation}
and so on.  If $\bfF(z)$ is analytic for $D_r(0)$ it follows that $\bfG(z) = \bfg(0) + \bfg(1)z + \bfg(2)z^2 + \cdots$ is also analytic for $z \in D_r(0)$.  Solving (\ref{e:(C.2)}) gives
\begin{eqnarray}
\label{e:(C.3)}
\lefteqn{\bff(0) = \bfg(0), \quad\bff(1) = \bfg(1) + [1 - p/(q+1)]\bfg(0), } \hspace{1cm} \nonumber \\
& & \bff(2) = \bfg(2) + [1 - p/(q+2)]\bfg(1) + [1 - p/(q+2)][1 - p/(q+1)]\bfg(0), \cdots \hspace{1cm}
\end{eqnarray}
and so on.  If $\bfG(z)$ is analytic for $z \in D_r(0)$ then $\bfF(z) = \bff(0) + \bff(1)z + \bff(2)z^2 + \cdots$ is also analytic for $z \in D_r(0)$.  In general we have
\begin{eqnarray*}
\bff(t) & = & \bfg(t) + [1 - p/(q+t)]\bfg(t-1) + \cdots + [1 - p/(q+t)] \cdots [1 - p/(q+1)]\bfg(0) \\
& = & \bfg(t) + \mbox{$ \sum_{s=1}^t$} [1 - p/(q+t)] \cdots [1 - p/(q+t-s+2)] [1 - p/(q+t-s+1)] \bfg(t-s)
\end{eqnarray*}
for all $t \in {\mathbb N}-1$.  Thus we now have a general formula for the desired inverse operator.  To justify (\ref{e:(C.1)}) we argue as follows. If $z \in D_1(0)$ and we write $(1-z)^{-1} = 1 + z + z^2 + \cdots$ then for each $k \in {\mathbb N}-1$ direct calculation gives $(1 - z)^{-1} [z^k] = \mbox{$\sum_{s \in {\mathbb N}-1}$} z^{s+k}$ from which it follows that
$$
[(1-z)^{-1}W](1-z)^{-1}[z^k] = p \cdot \mbox{$\sum_{s \in {\mathbb N}-1}$} \left[ \mbox{$\sum_{1 \leq i_1 \leq s+1}$} \mbox{$\prod_{j=1}^1$} (q+k+i_j)^{-1} \right]z^{s+1+k},
$$
$$
[(1-z)^{-1}W]^2(1-z)^{-1} [z^k] = p^2 \cdot \mbox{$\sum_{s \in {\mathbb N}-1}$} \left[ \mbox{$\sum_{1 \leq i_1 < i_2 \leq s+2}$}\, \mbox{$\prod_{j=1}^2$} (q+k+i_j)^{-1} \right]z^{s+2+k},
$$
and so on. An elementary induction can now be used to show that
$$
[(1-z)^{-1}W]^n(1-z)^{-1}[z^k] = p^n \cdot \mbox{$\sum_{s \in {\mathbb N}-1}$} \left[ \mbox{$\sum_{1 \leq i_1 < \cdots < i_n \leq s+n}$}\, \mbox{$\prod_{j=1}^n$} (q+k+i_j)^{-1} \right] z^{s+n+k}
$$
for all $n \in {\mathbb N}$ provided $z \in D_1(0)$.  If we calculate $[(1-z)I + W]^{-1}[z^k]$ but extract only the coefficients of $z^t$ for some fixed value of $t \in {\mathbb N}-1$ we need only consider values of $k \leq t$.  Thus we have
$$
\left( (1-z)^{-1}[z^k] \right)_t = 1,
$$
$$
\left( (-1)[(1-z)^{-1}W](1-z)^{-1}[z^k] \right)_t = (-1) p \cdot \mbox{$\sum_{1 \leq i_1 \leq t-k}$}\, \mbox{$\prod_{j=1}^1$} (q + k + i_j)^{-1},
$$
$$
\left( (-1)^2[(1-z)^{-1}W]^2(1-z)^{-1}[z^k] \right)_t = (-1)^2 p^2 \cdot \mbox{$\sum_{1 \leq i_1 < i_2 \leq t-k}$}\, \mbox{$\prod_{j=1}^2$} (q + k + i_j)^{-1},
$$
$$
\left( (-1)^3[(1-z)^{-1}W]^3(1-z)^{-1}[z^k] \right)_t = (-1)^3 p^3 \cdot \mbox{$\sum_{1 \leq i_1 < i_2 < i_3 \leq t-k}$}\, \mbox{$\prod_{j=1}^3$} (q + k + i_j)^{-1},
$$
and so on. In general
$$
\left( (-1)^{t-k}[(1-z)^{-1}W]^{t-k} [z^k] \right)_t = (-1)^{t-k} p^{t-k} \cdot \mbox{$\sum_{1 \leq i_1 < \cdots < i_{t-k} \leq t-k}$}\, \mbox{$\prod_{j=1}^{t-k}$} (q + k + i_j)^{-1}.
$$
Finally $\left((-1)^n[(1-z)^{-1}W]^n [z^k] \right)_t = 0$ for $n > t-k$. Collecting all of the terms shows that
\begin{eqnarray*}
\left( S[z^k] \right)_t & = & [1 - p/(q+k+1)][1 - p/(q+k+2)] \cdots [1 - p/(q+t)] \\
& = & \mbox{$\prod_{s=1}^{t-k}$} [1 - p/(q + k + s)]
\end{eqnarray*}
for each $t,k \in {\mathbb N}$ with $k \leq t$.  It follows that
\begin{eqnarray*}
\lefteqn{ S[ \bfg(0) + \bfg(1)z + \bfg(2)z^2 + \cdots + \bfg(t)z^t ] } \\
& \boldsymbol{\rightarrow} & \left\{ \bfg(t) + \mbox{$ \sum_{s=1}^t$} [1 - p/(q+t)] \cdots [1 - p/(q+t-s+2)] [1 - p/(q+t-s+1)] \bfg(t-s) \right\} \\
& = & \bff(t)
\end{eqnarray*}
for all $t \in {\mathbb N}-1$.  Therefore $S\bfG(z) = \bfF(z)$. It follows that $S \equiv [(1-z)I + W]^{-1}$ for $z \in D_1(0)$.  If $\bfF(z)$ is analytic for $z \in D_{1 + \epsilon}(0)$ where $\epsilon > 0$ the definition of $W\bfF(z)$ as an analytic function remains valid on this region. The difficulty with the definition of $S = [(1-z)I + W]^{-1}$ arises because if we define $\bfH(z) = \bfF(z) z^q$ and write
$$
H(z) = \bfH(1) + \bfH^{\, \prime}(1)(z-1) + \cdots
$$
for $z \in D_{1 + \epsilon}(0)$ then
$$
W(1 - z)^{-1}\bfF(z) = (p/z^q) \mbox{$\int_0^z$} \left[ \bfH(1)/(1 - \zeta) + \cdots \right] d\zeta = (-1)(p/z^q) \bfH(1) \log_e (1-z) + \cdots
$$
contains the term $\log_e (1 - z)$ which can only be defined as a single-valued analytic function if a suitable cut terminating at $z=1$ is inserted into the complex plane.  For instance the cut $\{ z \in {\mathbb C} \mid \Re(z) \in [1, \infty), \Im(z) = 0 \}$ is suitable.  Our earlier intuitive analysis suggested that $[(1-z)I + W]^{-1}$ has an {\em essential} singularity at $z=1$.  We can now see that it is not an {\em isolated} singularity.

\end{document}